\documentclass[final,10pt]{elsarticle}

\usepackage{amsmath}
\usepackage{amssymb}
\usepackage{amsthm}
\usepackage{amscd}
\usepackage{amsfonts}

\usepackage{graphicx}
\usepackage{color}
\usepackage{psfrag} 
\usepackage{comment,float,url}

\usepackage{amsmath}
\usepackage{soul,ulem}
\setstcolor{red}

\usepackage[top=2.8cm,bottom=2.8cm,left=2.5cm,right=2.5cm]{geometry}
\usepackage[default]{cantarell} %% Use option "defaultsans" to use cantarell as sans serif only
\usepackage[T1]{fontenc}

\newcommand{\bbb}{\textbf{b}}

\newcommand{\bff}{\textbf{f}}
\newcommand{\bg}{\textbf{g}}
\newcommand{\bh}{\textbf{h}}

\newcommand{\bn}{\textbf{n}}

\newcommand{\bu}{\textbf{u}}
\newcommand{\bv}{\textbf{v}}
\newcommand{\bw}{\textbf{w}}
\newcommand{\bx}{\textbf{x}}

\newcommand{\bphi}{\boldsymbol{\phi}}

\newtheorem{lem}{Lemma}

\newtheorem{rem}{Remark}
\newcommand{\pa}{\partial}
\newcommand{\f}{\frac}

\newcommand{\stress}{\mathcal{T}}

\usepackage[dvipsnames]{xcolor}
\usepackage{tikz}

\usepackage{newverbs}

\DeclareMathOperator*{\argmin}{arg\,min}

\begin{document}

\begin{frontmatter}
\title{El-WaveHoltz: A Time-Domain Iterative Solver for Time-Harmonic Elastic Waves}
\author[auth1]{Daniel Appel\"{o}}
\ead{appeloda@msu.edu}

\author[auth2]{Fortino Garcia\corref{corauthor}}
\ead{Fortino.Garcia@colorado.edu}

\author[auth2]{Allen Alvarez Loya} 
\ead{Allen.AlvarezLoya@colorado.edu}
\author[auth3]{Olof Runborg}
\ead{olofr@kth.se}
\address[auth1]{Department of Computational Mathematics, Science, and Engineering; Department of Mathematics, Michigan State University, East Lansing, MI 48824, USA}
\address[auth2]{Department of Applied Mathematics, University of Colorado, Boulder, CO 80309, USA}
\address[auth3]{Department of Mathematics, KTH, 100 44, Stockholm, Sweden}

\cortext[corauthor]{Corresponding author}

\begin{abstract}
We consider the application of the WaveHoltz iteration 
to time-harmonic elastic wave equations with energy conserving boundary conditions.
The original WaveHoltz iteration for acoustic Helmholtz problems is a fixed-point iteration
that filters the solution of the wave equation with time-harmonic forcing and boundary data.
As in the original WaveHoltz method,
we reformulate the fixed point iteration as a positive definite linear system of equations 
that is iteratively solved by a Krylov method.
We present two time-stepping schemes, one explicit and 
one (novel) implicit, which \textit{completely remove} time discretization error from the 
WaveHoltz solution by performing a simple modification of the initial data and time-stepping
scheme. Numerical experiments indicate an iteration scaling similar to that of the original WaveHoltz
method, and that the convergence rate is dictated by the shortest (shear) wave speed of the 
problem. We additionally show that the implicit scheme can be advantageous in practice for meshes with
disparate element sizes.
\end{abstract}

\end{frontmatter}

\section{Introduction}
Time-harmonic wave propagation problems are notoriously difficult to solve
by direct or  iterative methods due to the resolution requirements and the indefinite 
nature of the differential operator, especially at high frequencies.
For applications such as in solid mechanics, seismology and geophysics, 
we consider the time-harmonic elastic wave equation (or Navier 
equation \footnote{This is sometimes also referred to as the Navier-Cauchy 
equation, which is not to be confused with the ubiquitous Navier-Stokes equation.})
\begin{equation*}
\rho\omega^2 \bv + \nabla \cdot \mathcal{T}(\bv) = \bff(\bx), \ \  \bx \in \Omega,
\end{equation*}
for a domain $\Omega$ and frequency $\omega$. Here $\rho$ is the density,
$\bv \in \mathbb{R}^d$ is the displacement 
vector, $d$ the spatial dimension, $\mathcal{T}$ is the stress tensor,
and $\bff$ is the forcing. Few effective solvers and preconditioners are available for the Navier equation,
which are generally extensions of methods originally designed for 
the acoustic time-harmonic wave equation more commonly known as the 
Helmholtz equation. The efficient solution of the Helmholtz equation
via iterative methods is an active area of research with a variety of 
methods in both the frequency and time-domain. 
We refer to our previous paper \cite{WaveHoltz} 
for a more in-depth overview of the literature on techniques for 
solving the Helmholtz equation, as well as the 
review articles \cite{ernst2012difficult,Gander_Zhang_SIAM_REV,erlangga2008advances}.

The elastic wave equation models both pressure and shear waves and, 
as is the case for the Helmholtz equation,
the system of equations results in a discretization that is
highly indefinite for large frequencies. As for any wave propagation problem the resolution must increase 
with the frequency, and here the most stringent resolution constraint comes from the (shorter) shear wave wavelength. This, in tandem with
$d$ times the number of unknowns leading
to larger storage requirements, necessitates parallel, memory lean,
and scalable solvers that must be high order accurate
to mitigate dispersive errors, \cite{KreOli72}, causing the so-called \textit{pollution effect} \cite{pollution_error}.

While most methods have traditionally focused on solving the Helmholtz equation
in the frequency domain (we provide a review of some of these below), an alternative approach is to instead construct
iterative solvers in the time-domain. One such method, the 
so-called Controllability Method (CM), was first proposed by
Bristeau et al. \cite{bristeau1998controllability} and has recently received renewed interest in a series of papers by Grote et al.\cite{grote2019controllability,GROTE2020112846,Grote2021controllability}. The CM was extended 
to elastic media in \cite{monkola2008time,hoe2021solving}. The unknown in the CM is the initial data to the wave equation. In the CM this initial data is adjusted so that it produces an approximation to the  Helmholtz equation by solving a constrained
least-squares minimization where the objective function measures the 
deviation from time-periodicity. The minimization can be efficiently implemented using the conjugate gradient method, where the gradient is computed by solving the adjoint wave equation backwards in time. 

It is possible, however, to find a time-harmonic wave equation solution by solving a single wave equation
forward in time in each iteration. This can be done via the WaveHoltz method introduced 
in \cite{WaveHoltz} for the scalar wave equation. We now provide an overview of the WaveHoltz method which
was directly inspired by recent work on the CM, \cite{grote2019controllability}.

\subsection{Overview of the El WaveHoltz Iteration}

As in the CM, the WaveHoltz iteration (and the Elastic version we denote El WaveHoltz) iteratively updates the initial data to the wave equation. The marked difference between the two methods is that the WaveHoltz iteration updates the initial data by filtering the wave equation solution over one period (or an integer number of periods). The filtered solution is then used as the next initial data and thus the WaveHoltz method only requires one wave solve per iteration while the CM requires two. 

In \cite{WaveHoltz} we show that the (linear) iteration is convergent in both the continuous and discretized setting and that, if formulated as a linear system of equations, the underlying matrix is positive definite. We also showed that with energy conserving boundary conditions (Dirichlet or Neumann) the matrix is symmetric as long as the numerical method is symmetric (or symmetrizable). 

We emphasize that the filter used in the WaveHoltz method corresponds to a bounded linear operator with an inverse that is also bounded. Therefore, the number of iterations (and the condition number of the problem) is essentially independent of the gridsize $h$ for well resolved grids. This is in contrast to methods that discretize and solve the PDE directly. Such methods typically have a condition number that scales as $h^{-2}$ which makes it increasingly difficult to solve the problem as the solution becomes more accurate.       

The analysis in \cite{WaveHoltz} predicted that the WaveHoltz method in $d$ dimensions and accelerated by the conjugate gradient method converges to a fixed tolerance in $\mathcal{O}(\omega^d)$ iterations for energy conserving problems and numerical experiments indicated that it converges in $\mathcal{O}(\omega)$ iterations for open problems. The analytical predictions from \cite{WaveHoltz} are expected to hold here as well and in the experiments we carry out below we observe $\mathcal{O}(\omega^d)$. All these results and observations are independent of grid resolution indicating that our method can be particularly suitable when accurate solutions are required.

In this paper we focus solely on energy conserving boundary conditions (Dirichlet or normal stress) and leave the cases of impedance and non-reflecting boundary conditions to future work. In addition to introducing El WaveHoltz, we present several new results that are also retroactively applicable to our earlier work on WaveHoltz for the scalar wave equation \cite{WaveHoltz} and Maxwell's equations \cite{peng2021waveholtz}. 

Following the ideas of Stolk \cite{stolk2020timedomain} we introduce two new two-level time-stepping schemes -- one explicit and one implicit -- that remove the time-stepping error from the WaveHoltz solution. When either of these time-stepping methods are used the solution to the discrete WaveHoltz method is identical to the solution obtained by directly discretizing the frequency domain equation.   

For high frequency, large scale problems, parallel solution of the Navier equation is the only feasible option. For a parallel solver to scale well the ratio of communication to computation should be small. In general, there are two types of communication: a) the local communications between processors to exchange local degrees of freedom needed for stencil operations in the discretization of spatial derivatives, and b) global all-to-all operations such as computing the inner product between two global vectors. The WaveHoltz method has an intrinsic advantage compared to methods that work directly with the frequency domain equation in that the all-to-all communication that is required to update search directions in CG, GMRES etc. only needs to be computed once per $T = 2\pi/\omega$-period. Here we explore the effect of filtering over an additional number of periods to further reduce the number of all-to-all communications.         

We believe that the method we propose here is an attractive alternative to previously proposed methods.
In particular, El WaveHoltz is easily implemented if an elastic wave equation solver is already available. 
As we show in the numerical experiments section, El WaveHoltz can be one to two 
orders of magnitude faster compared to an algebraic multigrid (AMG) preconditioned 
GMRES solver for the frequency domain equation when using the symmetric interior 
penalty discontinuous Galerkin implementation available in MFEM~\cite{mfem}. 
There are, of course, many other solvers available;
the question of which method will be most efficient will 
depend on the details of the problem to be solved.
We now review some of the methods available in the literature.

\subsection{Literature Review}
One of the most common preconditioners for acoustic problems is the 
shifted Laplacian preconditioner (SLP), a more thorough review of which 
can be found in the review article by Erlangga \cite{erlangga2008advances}.
One of the first extensions of the SLP
to elastic media was introduced by Airaksinen et al.
\cite{AiraksinenDamped}, in which a finite element spatial discretization
for the preconditioned system is inverted by AMG. 
A more traditional finite difference multigrid SLP with line-relaxations 
was considered by Rizzuti and Mulder~\cite{RIZZUTI201647}. 
For both of these previous approaches, the effectiveness of a straightforward 
SLP is degraded for nearly incompressible media 
as the prolongation operators struggle to approximate the 
nullspace of the grad-div operator. To address this, a more recent extension 
was done by Treister \cite{treister2018shifted} in which a mixed-formulation
of the Navier equation is considered. 
While nearly incompressible media could be handled by the methods of
\cite{treister2018shifted}, this comes at the cost of doubling the number of
unknowns as well as additional storage requirements for precomputing the 
inverse of relaxation operators.

Another important class of methods for the solution of the 
Helmholtz equation are domain decomposition (DD) methods, for which
we refer the reader to \cite{Gander_Zhang_SIAM_REV} for a review.
In the short article \cite{BrunetSchwarzElastic}, it was shown 
that a classic Schwarz DD with overlap for elastic problems
converges for high frequencies, diverges for medium frequencies, and 
stagnates for small frequencies. Moreover, overlapping DD as a 
preconditioner for a GMRES accelerated solver exhibits convergence
behavior that depends strongly on the frequency $\omega$ with 
degrading performance for increasing frequency. To remedy this,
Brunet et al. introduced more general transmission conditions at the 
boundaries of overlapping domains in \cite{brunet2020natural}. These
transmission conditions, together with a 
sufficiently large enough overlap, yield convergence of the DD method
for all frequencies with the exception of $\{\omega/C_{\rm s}, \omega/C_{\rm p}\}$,
where $C_{\rm s}$ and $C_{\rm p}$ are the shear and pressure wave speeds, respectively.

For unbounded problems one of the most promising classes of preconditioners for the Helmholtz equation are the so-called sweeping preconditioners
by Engquist and Ying \cite{engquist2011sweepingH,engquist2011sweepingPML}.
These preconditioners construct an $LDL^T$ decomposition by sweeping through  
the domain layer-by-layer, with the key observation that the application of 
the Schur complement matrices found in the block diagonal matrix $D$ 
is equivalent to solving a quasi-1D(2D) problem in 2D(3D). In contrast to the
acoustic case, however, the sweeping preconditioner for 
time-harmonic elastic waves, \cite{tsuji2014sweeping}, exhibited an increase in the 
number of iterations with frequency for a heterogeneous 
media as the moving perfectly matched layer (PML) does not approximate 
Green's function as well. We note that the stable construction of PML 
for many elastic problems is still considered an open question 
\cite{BecFauJol03,AppKre05}. Similar to the sweeping preconditioner,
Belonosov et al. \cite{10.1190/geo2017-0710.1}
construct a preconditioner in 3D with damping that sweeps through
the domain along a coordinate axis while additionally 
homogenizing the medium in each layer. The preconditioner of \cite{10.1190/geo2017-0710.1}  
is inverted using FFT's and is accelerated with BiCGSTAB in the outer loop. 
As with the sweeping preconditioner, the choice of 
sweeping direction is important. Thus for problems where heterogeneity
is present in all directions this preconditioner is less effective.
Yet another solver with a sweeping nature is an extension of the 
Gordon and Gordon \cite{gordon2013robust} CARP-CG method
for Helmholtz problems to elastic media \cite{LiCarpCG}.
Despite its simplicity this method requires a large number 
of iterations, especially for heterogeneous media or 
problems with higher Poisson ratios. It should be emphasized that, 
although successful for unbounded problems, the efficiency of sweeping 
methods for energy conserving boundary conditions has largely not been 
demonstrated and their parallel implementation remains cumbersome.     

Instead of the $LDL^T$ decomposition used by the sweeping preconditioner,
other approaches constructing LU/ILU factorizations and preconditioners
are available. In \cite{el2011wavelet} an ILU preconditioner 
based on wavelet transforms with Gibbs reordering is used in a 
GMRES accelerated solver (with restarts) for time-harmonic elastic waves.
Wang et al. introduced a structured multifrontal 
algorithm using nested dissection based domain decomposition, together with 
hierarchical semi-separable (HSS)
compression for frontal matrices with low off-diagonal ranks in
\cite{wang2012massively}. The use of multilevel sequentially semi-separable 
(MSSS) matrix structure of the discretized elastic wave equation on 
Cartesian grids was leveraged in \cite{baumann2018fast} inside of 
an induced dimension reduction (IDR) accelerated ILU preconditioner.
The drawback of LU/ILU methods for the Navier equation is
the growth in memory and storage requirements.

The rest of this paper is organized as follows. In Section~\ref{sec::GovEq} we present
the Navier and elastic wave equations. In Section~\ref{sec::ElWHIter} we introduce the
WaveHoltz iteration applied to elastic problems with Dirichlet and/or 
free surface boundary conditions. In Section~\ref{sec::NumMeth_DiscAnalysis} we outline the numerical methods used to solve the elastic wave equation and present new results on time-stepping and Krylov acceleration. Numerical 
examples are presented in  Section~\ref{sec::NumericalExperiments}. Finally, we summarize and conclude in Section~\ref{sec::Conclusion}.
\section{Governing Equations}\label{sec::GovEq}

\subsection{The Elastic Wave Equation}\label{sec:elastic-wave-eqation-time}
The linear elastic wave equation in an isentropic material described by the density $\rho(\bx,t)$, the Lam\'{e} parameters $\mu(\bx)>0$ and $\lambda(\bx)>0$, and with a time-harmonic forcing takes the form 
\begin{equation} \label{eq:elastic-time}
\rho \bu_{tt} = \nabla \cdot \stress(\bu) - \text{Re}\{\bff(\bx) e^{i \omega t}\}, \ \  \bx \in \Omega, \ \ 0 \le t \le T. 
\end{equation}
Here $\bu = (u(\bx,t), v(\bx,t), w(\bx,t))$ is the displacement vector, $\bx = (x,y,z)^T$ is the Cartesian coordinate and $t$ is time. The stress tensor $\stress(\bu)$ can be decomposed into 
\begin{equation}
\stress(\bu) = \lambda (\nabla \cdot \bu) I + 2 \mu \mathcal{D}(\bu),
\end{equation}
where $\mathcal{D}(\bu)$ is the symmetric part of the displacement gradient 
\begin{equation}
\mathcal{D}(\bu) = \f{1}{2} \left(
\begin{array}{ccc}
2u_x & u_y+v_x & u_z+w_x\\
u_y+v_x & 2 v_y & v_z + w_y\\
u_z+w_x &  v_z + w_y & 2 w_z
\end{array}
\right).
\end{equation}

We assume that equation (\ref{eq:elastic-time}) is closed by time-harmonic boundary conditions specifying the displacement
\begin{equation}
\bu(\bx,t) = \text{Re}\{\bg(\bx) e^{i\omega t}\}, \ \ \bx \in \partial \Omega_{\rm D},
\end{equation}
or the normal stress 
\begin{equation}
\stress(\bu)\bn = \text{Re}\{\bh(\bx) e^{i\omega t}\}, \ \ \bx \in \partial \Omega_{\rm S},
\end{equation}
along with initial conditions 
\begin{equation}
\bu(\bx,0) = \bu_0(\bx), \ \ \frac{\partial \bu(\bx,0)}{\partial t} = \bu_1(\bx). 
\end{equation}

Multiplying (\ref{eq:elastic-time}) by $\bu^T$, integrating over $\Omega$ and invoking the divergence theorem yields  the energy estimate 
\begin{equation}
\frac{1}{2} \frac{d}{dt} \left( \| \sqrt(\rho) \bu_t \|^2 + \int_\Omega \lambda (\nabla \cdot \bu) I + 2 \mu (\mathcal{D}:\mathcal{D}) d \bx   \right) = - \int_\Omega   \cos (\omega t) \bu^T \bff(\bx) d \bx + \int_{\partial \Omega} \bu_t^T \stress(\bu) \bn \, dS.
\end{equation}
Here $\bn$ is the outward unit normal and the notation $(\mathcal{A}:\mathcal{B}) = \sum_{i = 1}^d \sum_{j=1}^d a_{i,j} b_{i,j}$ is the standard tensor contraction over two indices.

Thus, when there is no forcing, $\bff(\bx)=0$, the energy is conserved in time as long as $\bu_t^T \stress(\bu) \bn = 0$ on the boundary  $\partial \Omega$. The condition $ \stress(\bu) \bn = 0$ indicates that the boundary is stress free or free of traction. The Dirichlet condition on the velocity $\bu_t= 0$ also holds if the displacement vanishes for all time on the boundary, i.e. $\bu= 0$.

\begin{rem}
In the rest of this paper, unless otherwise noted, we will assume that the equations have been non-dimensionalized and that $\rho = 1$. 
\end{rem}

\subsection{The Time-Harmonic Elastic Wave Equation}\label{sec:elastic-wave-eqation-frequency}
Note that if the initial data of the elastic wave equation \eqref{eq:elastic-time} gives rise to a solution of the form $\bu(\bx,t) = \text{Re}\{\bv(\bx) e^{i\omega t}\}$ then $\bv$ satisfies the frequency domain equation
\begin{equation} \label{eq:elastic-frequency}
\omega^2 \bv + \nabla \cdot \stress(\bv) = \bff(\bx), \ \  \bx \in \Omega, 
\end{equation}
with boundary conditions on either the displacement or normal stress
\begin{equation}
\bv(\bx) = \bg(\bx), \ \ \bx \in \partial \Omega_{\rm D}, \quad \stress(\bv)\bn = \bh(\bx), \ \ \bx \in \partial \Omega_{\rm S},
\end{equation}
where $\stress$ is the stress tensor (and defined in the previous section).

For notational convenience we will refer to this as the elastic Helmholtz or, when there is no ambiguity, simply the Helmholtz equation though it is often called the Navier or Navier-Cauchy equation. We note that, in general, the Helmholtz solution $\bv$ is complex-valued. However, for boundary conditions that conserve the energy (such as Dirichlet and conditions on the normal stress) the corresponding solution $\bv$ becomes real-valued. For real-valued solutions, the corresponding time-harmonic solution of the elastic wave equation \eqref{eq:elastic-time} then simplifies to \mbox{$\bu(\bx,t) = \bv(\bx)\cos(\omega t)$}. The El-WaveHoltz method can be used to find the solution $\bv$ in both cases, but as we exclusively consider the energy conserving case here we describe the method for that case.  

%%%%%%%%%%%%%%
%
%%%%%%%%%%%%%%
\section{The El-WaveHoltz Iteration}\label{sec::ElWHIter}
The El-WaveHoltz iteration is a direct generalization of the WaveHoltz iteration introduced and analyzed in \cite{WaveHoltz}. In each iteration the elastic wave equation 
(\ref{eq:elastic-time}) is solved over one period $t\in[0,T]$, where $T=2\pi/\omega$.
Precisely, 
if we consider the energy conserving case, applying the WaveHoltz operator component wise to the initial displacement vector $\bu_0$ defines the El-WaveHoltz operator 
\begin{equation} \label{eq:el_WHI_OP}
 \Pi 
 \bu_0\\
 = 
 \frac{2}{T} \int \limits_{0}^{T} \left(\cos(\omega t) - \frac{1}{4} \right)
\bu(x,t)\\
 \ dt.
\end{equation}  
Here %$T = \frac{2\pi}{\omega}$ and 
$\bu(\bx,t)$ is the solution to (\ref{eq:elastic-time}) with the initial data 
$\bu(\bx,0)=\bu_0$ and $\frac{\partial\bu(\bx,0)}{\partial t} = 0$.
(For the energy conserving case one always uses $\bu_1= 0$). 
To briefly motivate \eqref{eq:el_WHI_OP}, we mention two consequences
of the particular form of the above operator. First, the filter \eqref{eq:el_WHI_OP}  is designed such that 
solutions of the elastic wave equation of the form \mbox{$\bu(\bx,t) = \bv(\bx)\cos(\omega t)$}
yield $\bv(\bx)$ as a fixed point. Second, it can be shown (see \cite{WaveHoltz}) that the constant shift of $-1/4$
guarantees that the operator \eqref{eq:el_WHI_OP} has a \textit{unique} fixed point.
The fixed point iteration then proceeds as
$$
  \bu^{(i+1)}=\Pi\bu^{(i)}, \qquad \bu^{(0)}=0,
$$
where we use $\bu^{(i)}$ to denote the $i^\text{th}$ WaveHoltz iterate.

As the analysis of this operator is the same as that for the scalar operator analyzed in \cite{WaveHoltz}, we will not repeat the analysis in detail here. Instead, we now highlight its most important properties. The first thing to note is that 
the elastic Helmholtz solution $\bv(\bx)$ is a fixed point of the operator.
Indeed,
if \mbox{$\bu(\bx,t) = \bv(\bx) \cos (\omega t)$} (and thus $\bu_0(\bx) = \bv(\bx)$), then the integral in (\ref{eq:el_WHI_OP}) can trivially be evaluated 
\begin{equation} 
\Pi \bv(x) =  \frac{2}{T} \int \limits_{0}^{T} \left(\cos(\omega t) - \frac{1}{4} \right)
\cos(\omega t)\bv(x)\\
 \ dt = \bv(x),
\end{equation}  
%showing that the elastic Helmholtz solution is a fixed point of the operator. 
Further, we denote by $\mathcal{S}$ the operator $\Pi$ for the case when $\bff = 0$. If $(\lambda_j^2, \bphi_j)$ is the eigendecomposition satisfying $ \lambda_j^2 \bphi_j =  \nabla \cdot \mathcal{T}(\bphi_j)$, then for a general initial displacement the solution will be of the form $\sum_{j=0}^{\infty} d_j \cos(\lambda_j t) \bphi_j$. Defining
\[
 \beta(\lambda) \equiv \frac{2}{T} \int \limits_{0}^{T} \left(\cos(\omega t) - \frac{1}{4} \right)
\cos(\lambda t) dt,
\]
we obtain as in \cite{WaveHoltz} that $\mathcal{S}$ can be expressed as
\[
\mathcal{S} \sum_{j=0}^{\infty} d_j  \bphi_j \equiv \sum_{j=0}^{\infty} \beta(\lambda_j) d_j  \bphi_j.
\]
%which gives the filtered solution to the elastic wave equation when $\bff = 0$. 
If $\omega \neq \lambda_j$ for all $j$ then the spectral radius of $\mathcal{S}$ is given by $ \max_j | \beta(\lambda_j) | < 1$ (see Lemma 2.1 in \cite{WaveHoltz}) so the iteration will converge. Since the operator $\Pi$ is affine, we may find the fixed point (or equivalently the elastic Helmholtz solution) by solving the equation $(\mathcal{I} - \mathcal{S}) \bv \equiv \mathcal{A} \bv= \bbb \equiv \Pi {\bf 0}$. As is the case for the scalar Helmholtz equation, the eigenvalues of $\mathcal{A}$ lie in $(0,3/2)$ and the condition number scales with the frequency as ${\rm cond}(\mathcal{A}) \sim \omega^{2d}$ in $d$ dimensions. 

We emphasize that here $\mathcal{A}$ is a self-adjoint, positive definite and bounded operator. Thus once $\mathcal{A}$ is discretized it will be possible to apply the conjugate gradient method.
We can compute the action of $\mathcal{A}$ from the action of $\mathcal{S}$
which is obtained by solving the elastic wave equation. We
do not need to explicitly form $\mathcal{A}$.
 Moreover, as the condition number {\it does not} depend on the discretization size, the number of iterations are not expected to increase as the solution becomes more accurate due to grid refinement. We also note that since ${\rm cond}(\mathcal{A}) \sim \omega^{2d}$ the conjugate gradient method is expected to converge to a fixed tolerance in $\omega^d$ iterations.  

Finally, as mentioned above it is possible to define the iteration as the integral over multiple periods in order to reduce the number of all-to-all communication in the Krylov iteration. For example, if the number of periods is $K$ then we can define the filtering as  
 \begin{equation}
 \Pi_K \label{eq:many_periods_filter}
 \mathbf{u}_0\\
 = 
 \frac{2}{KT} \int \limits_{0}^{KT} \left(\cos(\omega t) - \frac{1}{4} \right)
\mathbf{u}\\
 \ dt,\  T = \frac{2\pi}{\omega}.
 \end{equation}

\begin{rem}
For general boundary conditions (e.g. non-reflecting or impedance), 
the iteration converges much faster, typically in $\sim \omega$ iterations independent of dimension.   
For this case
$\frac{\partial \bu(\bx,0)}{\partial t} = \bu_1(\bx)$ will not be zero and we must seek the initial data $\bu_0$ and $\bu_1$ simultaneously. The El-WaveHoltz operator then is
 \[
 \Pi 
 \begin{bmatrix}
 \bu_0\\
 \bu_1
 \end{bmatrix}
 = 
 \frac{2}{T} \int \limits_{0}^{T} \left(\cos(\omega t) - \frac{1}{4} \right)
 \begin{bmatrix}
\bu\\
\bu_t
 \end{bmatrix}
 \ dt,\  T = \frac{2\pi}{\omega}.
\]
This operator is more difficult to analyze; see \cite{WaveHoltz2}.
\end{rem}

\section{Numerical Methods and Discretization}\label{sec::NumMeth_DiscAnalysis}
An attractive feature of El-WaveHoltz is that it can be used together with any convergent discretization of the elastic wave equation. Here we consider the conservative curvilinear finite difference method from \cite{AppPet09} and the symmetric interior penalty discontinuous Galerkin method \cite{De-Basabe:2008aa,GSSwave}. We give a very brief description of these methods below and refer the reader to \cite{AppPet09,De-Basabe:2008aa} for details. 

Although highly non-intrusive, the one additional discretizational detail required by El-WaveHoltz is how to discretize the integral in (\ref{eq:el_WHI_OP}). As the integrand is periodic (once converged) we always use the trapezoidal rule.

\subsection{El-WaveHoltz by Finite Differences}\label{sec::ElWHI_FD}
To discretize the elastic wave equation 
%(\ref{sec:elastic-wave-eqation-time}) 
(\ref{eq:elastic-time})
in a general non-Cartesian geometry we write 
%(\ref{sec:elastic-wave-eqation-time}) 
(\ref{eq:elastic-time})
in a curvilinear coordinate system that conforms with the boundaries of the domain but that can be mapped back to the unit square (cube). Thus, we assume that there is a one-to-one mapping 
\begin{gather*}
 x=x(q,r),\ \ y= y(q,r), \ \ (q,r)\in [0,1]^2,
\end{gather*}
from the unit square to the domain of interest. Then the two dimensional version of (\ref{eq:elastic-time})
%(\ref{sec:elastic-wave-eqation-time}) 
becomes
\begin{multline*}
 J \rho \f{\pa^2 u}{\pa t^2} = 
 \f{\pa}{\pa q}\Big[ J q_x \big[ \left(2 \mu + \lambda \right)
 \left(q_x\pa_q+r_x\pa_r\right) u + \lambda \left(q_y\pa_q+r_y\pa_r\right) v \big] 
 + J q_y\big[ \mu \left(\left(q_x\pa_q+r_x\pa_r\right) v + \left(q_y\pa_q+r_y\pa_r\right) u  \right)
 \big] \Big] \\
 +\f{\pa}{\pa r}\Big[ J r_x \big[ \left(2 \mu + \lambda \right)
 \left(q_x\pa_q+r_x\pa_r\right) u + \lambda \left(q_y\pa_q+r_y\pa_r\right) v\big]
 + J r_y\big[ \mu \left(\left(q_x\pa_q+r_x\pa_r\right) v + \left(q_y\pa_q+r_y\pa_r\right) u  \right)
 \big] \Big],
\end{multline*}
\begin{multline*}
 J \rho \f{\pa^2 v}{\pa t^2} =
 \f{\pa}{\pa q}\Big[ J q_x \big[ 
 \mu \left( \left(q_x\pa_q+r_x\pa_r\right) v + \left(q_y\pa_q+r_y\pa_r\right) u  \right)
 \big]
 + J q_y\big[ 
  \left(2 \mu + \lambda \right)
     \left(q_y\pa_q+r_y\pa_r\right) v + \lambda \left(q_x\pa_q+r_x\pa_r\right) u
 \big] \Big] \\
 +\f{\pa}{\pa r}\Big[ J r_x \big[    
 \mu \left( \left(q_x\pa_q+r_x\pa_r\right) v + \left(q_y\pa_q+r_y\pa_r\right) u  \right)
 \big]  
 + J r_y\big[ 
 \left(2 \mu + \lambda \right)
 \left(q_y\pa_q+r_y\pa_r\right) v + \lambda \left(q_x\pa_q+r_x\pa_r\right) u
 \big] \Big].
\end{multline*}
Here $J=x_q y_r - x_r y_q$ is the Jacobian of the mapping. Also note that we have considered the case without forcing for brevity. 

We discretize the unit square $(q,r)\in [0,1]^2$ by a uniform grid on which we introduce real valued grid functions $[u_{i,j}(t),v_{i,j}(t)]=[u(q_i,r_j,t),v(q_i,r_j,t)]$. On this grid we apply an energy stable discretization 
\begin{equation}\label{eq:curve1disc}
  \rho J \f{\pa^2 u_h}{\pa t^2} =  L^{(u)}(u_h,v_h), \ \ \ \   \rho J \f{\pa^2 v_h}{\pa t^2} = L^{(v)}(u_h,v_h).
\end{equation}
Here $\rho J$ is a diagonal matrix containing the metric information and $u_h, v_h$ are vectors containing all the grid function values.
The (lengthy) exact definitions of $L^{(u)}(u_h,v_h) ,L^{(v)}(u_h,v_h)$ can be found in \cite{AppPet09}.   

To discretize the equations in time we either use the standard second order accurate centered differences, or one of the time-corrected schemes discussed below. For the standard second order accurate centered difference approximation in time, the fully discrete equations take the form
\begin{gather}\label{eq:4}
 \begin{aligned}
(\rho J)(u_h^{n+1}-2u_h^n+u_h^{n-1})= \Delta t^2 L^{(u)}(u_h^n,v_h^n),\\
(\rho J)(v_h^{n+1}-2v_h^n+v_h^{n-1})= \Delta t^2 L^{(v)}(u_h^n,v_h^n).
 \end{aligned}
\end{gather}
Then, if $(u,v)_{\rho J}$ is the weighted inner product defined by $(f,(\rho J)^{-1}\, g)_{\rho J}=(f,g)_h$, and  $C_e(t^{n+1})$ is the discrete energy
 \begin{equation}\label{eq:8}
   C_e(t^{n+1})= \| D^t_+ u^n\|_{\rho J}^2 + \| D^t_+ v^n\|_{\rho J}^2
   -( u^{n+1}, (\rho J)^{-1} L^{(u)}(u^n,v^n))_{\rho J}  -( v^{n+1}, (\rho J)^{-1} L^{(v)}(u^n,v^n))_{\rho J},
\end{equation}  
one can show that this discrete energy is conserved \cite{AppPet09}. 

Note that (\ref{eq:4}) is slightly non-symmetric and needs to be diagonally scaled to become symmetric. 
Here we scale by 2 along sides with free surface boundary conditions, and by 4 in corners where free surfaces meet. 
Incorporating this scaling through the multiplication by a scaling matrix $\Lambda$, the method can be formally written as 
\begin{equation}\label{eq:symmertic_FD}
M(\bu_h^{n+1} - 2\bu_h^{n} + \bu_h^{n-1}) = \Delta t^2 L_h \bu_h^n. 
\end{equation}
Here $M = {\rm diag}(\Lambda \rho J ,\Lambda \rho J)$ and $L_h = {\rm diag}(\Lambda L^{(u)},\Lambda L^{(v)} )$ are symmetric and $M$ is diagonal. However, as $M^{-1} L_h$ is not in general symmetric, the iteration (\ref{eq:el_WHI_OP}) will produce a symmetrizable but not symmetric operator. We will show below that this requires a minor modification of the conjugate gradient algorithm when used together with the iteration (\ref{eq:el_WHI_OP}).

\subsection{El-WaveHoltz by Symmetric Interior Penalty Discontinuous Galerkin Method}\label{sec::SIPDG}
As an alternative to the finite difference method outlined above, we will also consider the Symmetric Interior Penalty Discontinuous Galerkin (SIPDG) method \cite{De-Basabe:2008aa,GSSwave}. Let $\Omega_h$ be a finite element partition of the computational domain $\Omega$, with $\Gamma_h$ the set of all faces. Then \eqref{eq:elastic-time} can be reformulated 
into the interior-penalty weak formulation: Find $\bu_h \in (0,T) \times V_h$
such that
  \begin{align} \label{eq:SIPDG}
    \sum_{E \in \Omega_h} \left(\rho \frac{d^2 \bu_h}{dt^2}, \bv\right)_E + \sum_{E \in \Omega_h} B_E(\bu_h, \bv) + \sum_{\gamma \in \Gamma_h} J_\gamma(\bu_h, \bv; S, R) = - \cos(\omega t) \sum_{E \in \Omega_h} (\bff, \bv)_E,
  \end{align}
for all $\bv \in V_h$. Here 
  \begin{gather*}
    (\bu, \bv)_E  = \int_E \bu \cdot \bv \, dE, \\
    B_E(\bu, \bv) = \int_E \left[\lambda (\nabla \cdot \bu)(\nabla \cdot \bv) + \mu (\nabla\bu + \nabla\bu^T) : \nabla \bv\right] \, dE, \\
    J_\gamma(\bu, \bv; S, R) =   - \int_\gamma \{\stress(\bu) \bn \} \cdot [\bv] \, d\gamma 
                                  + S \int_\gamma \{\stress(\bv) \bn \} \cdot [\bu] \, d\gamma 
                                  + R \int_\gamma \{\lambda + 2\mu\} [\bu] \cdot [\bv] \, d\gamma,
\end{gather*}
% where $\nabla \bu : \nabla \bv = u_{i,j} v_{i,j}$,
where $\{\cdot \}$ and $[\cdot]$ denote the average and jump of a function,
respectively. The parameter $R$ is the penalty and $S$ determines 
the particular flavor of  IPDG. We %thus 
set $S=-1$, corresponding to the
Symmetric IPDG \cite{GSSwave}. In this case, $J_\gamma$ is 
symmetric with respect to $\bu_h$ and $\bv$ so that together 
with the symmetry of $B_E$ we have that the stiffness matrix is 
symmetric. Thus SIPDG provides a symmetric discretization
of the elastic wave equation, which will allow the use of conjugate
gradient to accelerate convergence of the El-WaveHoltz iteration.

Our solver is implemented in MFEM\footnote{\url{www.mfem.org}} \cite{mfem} and is essentially a direct extension of example 17 to the time-domain. 
Depending on the mesh, our choice of finite element space $V_h$ is typically one of two broken spaces. We choose either $\mathcal{P}^p(E)$, the space of polynomials of total degree at most $p$ on triangles, or $\mathcal{Q}^p(E)$,
the space of polynomials of at most degree $p$ on quadrilaterals. 
Unless otherwise noted, for the penalty parameter we make the choice $R = (p+1)(p+2)$
motivated by the analysis of \cite{shahbazi2005explicit}.

With the standard second order explicit time discretization, the matrix form of (\ref{eq:SIPDG}) becomes 
\[
M_\rho (\bu_h^{n+1} - 2 \bu_h^{n} +\bu_h^{n-1}) = \Delta t^2 \left[L_h \bu_h^{n}  - \bff_h \cos(\omega t_n)\right].  
\]
As for the finite difference method, $M_\rho^{-1} L_h$ is not (in general) symmetric and this will require a minor modification of the conjugate gradient algorithm when this scheme is used together with the iteration (\ref{eq:el_WHI_OP}).  

For this explicit time-stepping and the error corrected time-stepping discussed below, 
we use the ${\rm CFL}$ condition from \cite{el_dg_dath}
  \begin{align}\label{eqn::DG-CFL}
    \Delta t < \frac{{\rm CFL} \cdot h_{\rm min}}{(p + \frac{3}{2})^2 \sqrt{\frac{2\mu + \lambda}{\rho}}},
  \end{align}
where $h_{\rm min}$ is the smallest diameter of the elements and ${\rm CFL}$ depends on the time-stepper.
For the second order centered scheme, we typically choose ${\rm CFL}\sim$ 0.4--0.8.

\subsection{Explicit Time-Corrected Scheme}\label{sec::ExplicitTime-Stepping}
If the elastic Helmholtz equation (\ref{eq:elastic-frequency}) is discretized directly, the solution satisfies the equation (in this section we take $\rho = 1$ and for notational clarity we suppress the subscript $h$ for $\bv_h$ and $\bu_h$)  
\begin{equation} 
\omega^2 \bv + L_h \bv = \bff(\bx). 
\end{equation}

Consider the elastic wave equation time marched with the second order method
\begin{equation} \label{eq:elastic-time-explicit}
  \bu^{n+1} - 2\bu^n + \bu^{n-1} = \Delta t^2\left[ L_h \bu^n -\bff\cos (\omega t_n) \right],
\end{equation}
and started with the initial data
\begin{align*}
   \bu^0 = \bu_0,\qquad \bu^{-1} = \bu_0 - \frac{\Delta t^2}{2}L_h\left(\bu_0+\bff\right).
\end{align*}
Then, as we showed in \cite{WaveHoltz,WaveHoltz2} for the scalar wave equation, 
the fixed point iteration $\bu^{(i+1)} = \Pi \bu^{(i)}$ with $\bu^{(0)} = 0$
converges to $\bu^\infty$ which is a solution to
the elastic Helmholtz equation with a modified frequency
\begin{equation} 
\tilde{\omega}^2 \bu^\infty(\bx) + L_h \bu^\infty(\bx) = \bff(\bx),  \quad \tilde \omega = \frac{2\sin(\Delta t \omega/2)}{\Delta t}. 
\end{equation}
For this second order time discretization the difference between the final converged $\bu^\infty$ and the exact solution is $\mathcal{O}(\Delta t^2)$. Thus if a high order accurate spatial discretization is used, time discretization errors will limit the accuracy of the El-WaveHoltz solution. To reduce this error, a time discretization which is at least as accurate as the spatial discretization can be used. It is also possible, however, to use the technique proposed by Stolk in \cite{stolk2020timedomain} to modify the second order time-stepping method and eliminate the error altogether. The corrected scheme in \cite{stolk2020timedomain}, introduced as a time-domain preconditioner, is the straightforward modification 
\begin{equation} \label{eq:elastic-time-explicit-corrected}
   \bu^{n+1} - 2\bu^n + \bu^{n-1} = \frac{\tilde \omega^2}{\omega^2} \Delta t^2\left[L_h \bu^n -\bff\cos (\omega t_n)\right].
\end{equation}
As \cite{stolk2020timedomain} solves the equations in the frequency domain no initial data is used. Here, as we work in the time-domain, we must also modify the computation of $\bu^{-1}$ accordingly:
\begin{align*}
   \bu^0 = \bu_0,\qquad \bu^{-1} = \bu_0 - \frac{\tilde \omega^2}{\omega^2}\frac{\Delta t^2}{2}L_h\left(\bu_0+\bff\right).
\end{align*}

\subsection{Implicit Time-Corrected Scheme}\label{sec:whimp}
For a DG discretization, the use of an explicit time-stepping scheme
for the elastic wave equation requires a ${\rm CFL}$ condition that
shrinks as $\mathcal{O}(p^{-2})$ where $p$ is the polynomial order
within an element. For meshes with geometrical stiffness and DG discretizations of high order an implicit scheme can be used to avoid the restrictive time-step of an explicit scheme. 

To that end, we use the following second order, A-stable
implicit time-stepping scheme to solve the elastic wave
equation,
  \begin{align}\label{eqn::ImplicitScheme}
    \frac{\bu^{n+1}-2 \bu^n + \bu^{n-1}}{\Delta t^2} =
    \frac12 L_h(\bu^{n+1} + \bu^{n-1}) 
    -\bff \cos(\omega t_{n})\cos(\omega \Delta t).
  \end{align}
With a time-step $\Delta t=T/k$ for some integer $k$, the scheme \eqref{eqn::ImplicitScheme}
is completed by initial data 
$$
   \bu^0 = \bu_0,\qquad \bu^{-1} = \left(I - \frac{\Delta t^2}{2}L_h\right)^{-1}\left(\bu_0-\frac{\Delta t^2}{2}\bff\cos(\omega \Delta t)\right).
$$
As with the explicit method, El-WaveHoltz then converges to
$\bu^\infty$ which satisfies the elastic Helmholtz equation with a modified frequency
\begin{equation} 
\tilde{\omega}^2 \bu^\infty(\bx) + L_h \bu^\infty(\bx) = \bff(\bx),  \quad \tilde \omega = 2\frac{\sin(\omega \Delta t/2)}{\Delta t \sqrt{\cos(\omega \Delta t)}}.
\end{equation}
In a procedure similar to that done in the explicit case, it is possible 
to modify the implicit scheme to ensure the converged El-WaveHoltz 
solution is free from time discretization errors. This requires the 
modified scheme,
  \begin{align*}
    \frac{\bu^{n+1}-\alpha \bu^n + \bu^{n-1}}{\Delta t^2} =
    \frac12 L_h(\bu^{n+1} + \bu^{n-1}) 
    -\bff \cos(\omega t_{n})\cos(\omega \Delta t),
  \end{align*}
where 
  \begin{align}\label{eqn::alpha}
    \alpha = \cos(\omega \Delta t) (2 + \omega^2 \Delta t^2) \approx 2 - \frac{5 (\omega \Delta t)^4}{12} + \mathcal{O}(\Delta t^6),
  \end{align}
with modified initial data 
$$
   \bu^0 = \bu_0,\qquad \bu^{-1} = \left(I - \frac{\Delta t^2}{2}L_h\right)^{-1}\left(\frac{\alpha}{2}\bu_0-\frac{\Delta t^2}{2}\bff\cos(\omega \Delta t)\right).
$$

For the stability of the method it is necessary to have $|\alpha| < 2$. 
This choice of the time-step corresponds to a (mild) requirement of at least five 
time-steps per iteration (See details in \ref{appendix::RestrictionPlot}).

\begin{rem}\label{rem::Implicit_scheme}
The use of the above implicit scheme indeed allows one to circumvent a (potentially) restrictive ${\rm CFL}$
condition and take significantly larger time-steps compared to the explicit scheme. 
However, we note that a small number of time-steps can lead to inaccurate quadrature and values of the
discrete filter transfer function larger than one (see \ref{appendix::TrapQuad} for a
discussion and example). Thus if an eigenvalue of the discrete Laplacian is close
to resonance (i.e. close to $\omega$) and a small number of quadrature points are used,
then the linear system may become indefinite. In this case, we recommend using MINRES instead of 
CG to accelerate convergence of El-WaveHoltz. 

In the original WaveHoltz paper~\cite{WaveHoltz}, 
an application of Weyl asymptotics for the continuous problem eigenvalues gave that 
	\begin{align*}
		\delta \sim \omega^{-d},
	\end{align*}
where the relative gap to resonance, $\delta$, is such that 
	\begin{align*}
		\delta = \delta_{j^*}, \quad j^* = \argmin_j |\delta_j|, \quad \delta_j = \frac{|\omega - \lambda_j|}{\lambda_j},
	\end{align*}
and $\lambda_j$ is an eigenvalue of the Laplacian. Thus in higher dimensions there are more 
and more ``problematic'' modes close to resonance which would require a smaller time-step
for accuracy of the trapezoidal rule applied to the integral (\ref{eq:el_WHI_OP}). 
An option to mitigate the indefiniteness of the linear system due to 
an inaccurate trapezoidal rule would be to consider a higher order 
quadrature, or potentially identifying and removing problematic modes via deflation. 
\end{rem}

\begin{rem}
We  remark that it is also possible to remove the time discretization error by modifying the weights in the trapezoidal rule as in \cite{peng2021waveholtz}
\begin{align}
 \frac{2 \Delta t}{T}\sum_{n=0}^{N_{\rm t}} \frac{\cos(\omega t_n)}{\cos(\frac{2\sin(\Delta t \omega/2)}{\Delta t} t_n)} \left(\cos(\omega t_n)-\frac{1}{4}\right)  \bu^{n}.
\end{align}
It should however be noted that there is a risk that the denominator in this expression can become arbitrarily close to zero unless care is taken. We do not use this technique in any of the examples in this paper.
\end{rem}

\subsection{Krylov Solution of  the El-WaveHoltz Iteration}
Let $\Pi_h$ and ${\mathcal S}_h$ be the matrices corresponding to a discretization of the El-WaveHoltz operators $\Pi$ and ${\mathcal S}$ using either the finite difference or the SIPDG method. Then the iteration is (in this section a superscript $(i)$ denotes iteration and a superscript $n$ denotes time-step)
\begin{align*}
%\bu_h^0 &= \Pi_h {\bf 0},\\
\bu_h^{(0)} &= {\bf 0},\\
\bu_h^{(i+1)} &= \Pi_h \bu_h^{(i)}, \ \ i = 0,1,\ldots
\end{align*}
The solution to this fixed point iteration can also be found by solving 
\begin{equation} \label{eq:amat_non_sym}
(I-\mathcal{S}_h) \bu_h = \bbb \equiv \Pi_h {\bf 0}=\bu_h^{(1)},
\end{equation}
where the action of the matrix $(I-\mathcal{S}_h)$ requires (\ref{eq:elastic-time}) with $\bff = 0$ to be solved for one period, $T = 2\pi/\omega$, and the right hand side is pre computed by solving (\ref{eq:elastic-time}) with the $\bff$ at hand. 

Let ${\mathcal A}_h = I-\mathcal{S}_h$. We know from \cite{WaveHoltz} that the eigenvalues of ${\mathcal A}_h$ are in the interval $(0,3/2)$ so that ${\mathcal A}_h$ is positive definite. 
When applying ${\mathcal A}_h$ to $\bu$, we note that the methods for the elastic wave equation we consider here produce solutions $\{\bu^0,\bu^1,\ldots,\bu^{N_{\rm t}}\}$ at time instances $0,\Delta t, 2\Delta t, \ldots$, according to the recursion 
\begin{eqnarray} 
\bu^{-1} &=& a_0 \bu + a_1 M^{-1}W \bu, \label{eq:qrecursion1} \\
\bu^0 &=& \bu,  \label{eq:qrecursion2} \\
\bu^{n+1} &=& \kappa\bu^{n}-\bu^{n-1} + \gamma M^{-1}W \bu^n, \label{eq:qrecursion3}
\end{eqnarray}
for some diagonal matrix $M$, symmetric matrix $W$ and scalars $a_0$, $a_1$, $\kappa$, $\gamma$.
It follows that the matrix 
${\mathcal A}_h$ is in general not symmetric, as
\begin{equation} \label{eq:qoper}
{\mathcal A}_h \bu = \sum_{n=0}^{N_{\rm t}}
\alpha_n\bu^n={\mathcal P}(M^{-1}W)\bu,
\qquad \alpha_n = \begin{cases}
\frac{\Delta t}2, & \text{$n=0$ or $n=N_{\rm t}$}, \\
\Delta t, & \text{otherwise},
\end{cases}
\end{equation}
for some polynomial ${\mathcal P}$ of degree $N_{\rm t}-1$. 
%will not be symmetric even if $M$ and $S$ are. 
However, 
%as the operator $Q$ can be expressed as a polynomial $P_{Q}$ of degree $N_{\rm t}-1$ in $M^{-1}S$, 
$M{\mathcal A}_h$ is symmetric as each term in 
$M{\mathcal P}(M^{-1}W)$ is symmetric.
\begin{comment}
However, as the operator $Q$ can be expressed as a polynomial $P_{Q}$ of degree $N_{\rm t}-1$ in $M^{-1}S$, we have that  
\begin{equation}
Q \bw  = P_{Q}(M^{-1}S) \bw.
\end{equation}
As a result, the quantity $\bv^T MQ \bw$ can be written 
\begin{multline}
\bv^T MQ \bw = \bv^T M P_{Q}(M^{-1}S) \bw = \bv^T M (\beta_1 I + \beta_2 M^{-1}S + \ldots +  \beta_{N_{\rm t}} (M^{-1}S)^{N_{\rm t}-1}) \bw \\
= \beta_1 \bv^T M \bw  + \beta_2  \bv^T S \bw + \beta_3  (S \bv)^T M^{-1} (S \bw) \ldots +  \beta_{N_{\rm t}} \bv^T (M^{-1}S)^{N_{\rm t}-1}) \bw.
\end{multline}
The form of the last term in this expression depends on whether $N_{\rm t}$ is even or odd, and is either of the form 
$$
((M^{-1}S)^p\bv)^T S ((M^{-1}S)^p\bw),
$$
or 
$$
(S (M^{-1}S)^q\bv)^T M^{-1} (S(M^{-1}S)^q\bw).
$$
In either case, all terms in the expression for $\bv^T MQ \bw$ are symmetric with respect to $\bv$ and $\bw$. We formulate this observation as a lemma.  
\begin{lem}
The operator $Q$ defined by (\ref{eq:qoper}) and the recursion (\ref{eq:qrecursion1})--(\ref{eq:qrecursion3}) satisfies the identity
\begin{equation} \label{eq:same_scalar}
\bv^T MQ \bw = \bw^T MQ \bv.
\end{equation} 
\end{lem}

A consequence of this lemma is that rather than applying the conjugate gradient method to (\ref{eq:amat_non_sym}), we instead solve 
\end{comment}
Thus rather than applying the conjugate gradient method to (\ref{eq:amat_non_sym}), we instead solve 
\begin{align}\label{eq:amat_sym}
	M (I-{\mathcal S}_h) \bu_h = M\bbb.
\end{align}
We note that the main cost of applying the matrix ${\mathcal A}_h$ is in computing ${\mathcal S}_h\bu_h$. Since the matrix $M$ is diagonal for the finite difference method and block-diagonal for SIPDG, the difference in cost between applying \eqref{eq:amat_sym} over \eqref{eq:amat_non_sym} is small and amortized by the advantage of not having to store a Krylov subspace when using the conjugate gradient or conjugate residual method. 

\begin{rem}
In some of the experiments below we use the conjugate residual rather than conjugate gradient method. The reason for this is that it has the property that the residual is non-increasing, which we have found gives a predictable and robust iteration count when doing parameter sweeps over $\omega$. When conjugate gradient is used we sometimes observe that we get ``lucky'' and converge in very few iterations for a few frequencies. When considering practical applications it is of course good to have such luck, but as we are trying to present the average behavior of our method here we prefer conjugate residual.        
\end{rem} 

\section{Numerical Experiments}\label{sec::NumericalExperiments}
In this section we present numerical experiments that demonstrate the properties of the method.
We start with numerical experiments that demonstrate the spatial accuracy with and without the time-stepping correction for the finite difference and the discontinuous Galerkin method. Unless otherwise specified, the following computations were performed on Maneframe II at the Center for Scientific Computation at Southern Methodist University using a dual Intel Xeon E5-2695v4 2.1 GHz 18-core Broadwell processor with 45 MB of cache and 256 GB of DDR4-2400 memory.
\subsection{Accuracy of the Finite Difference Method}
We consider solving the elastic Helmholtz equation with Lam\'{e} parameters $\lambda = \mu = 1.0$, where the forcing function is chosen so that the displacements are given by   
\begin{align}\label{eqn::QuarticExample}
u = v = 16^2 x^2(x - 1)^2 y^2(y - 1)^2.
\end{align}
\begin{table}[ht]
\caption{$L_1, L_2$ and $L_{\infty}$ errors of the computed solution with corresponding estimated rates of convergence.\label{tab::ConvergenceFD}} 
\begin{center} 
 \begin{tabular}{|c|c|c|c|c|c|c|} 
\hline 
 n & $L_1$ error & Convergence & $L_2$ error & Convergence & $L_{\infty}$ error & Convergence \\ 
 \hline 
 20 & 3.86(-3) & - & 3.86(-3) & - & 3.86(-3) & - \\ 
40 & 9.21(-4) & 2.06 & 9.21(-4) & 2.06 & 9.21(-4) & 2.06 \\ 
80 & 2.25(-4) & 2.03 & 2.25(-4) & 2.03 & 2.25(-4) & 2.03 \\ 
160 & 5.55(-5) & 2.02 & 5.55e(-5) & 2.02 & 5.55(-5) & 2.02 \\
\hline 
\end{tabular} 
\end{center} 
\end{table} 

We take the frequency to be $\omega = 1.0$ and enforce Dirichlet boundary conditions on the boundary of the unit square $(x,y) \in [0,1]^2$. To verify accuracy, we set the tolerance to $10^{-15}$ for the conjugate residual method as the stopping criteria and compute the error in $u$ to the exact solution. We use the finite difference method of Section~\ref{sec::ElWHI_FD} together with the standard explicit second order time-stepping scheme, and verify the convergence of the method by grid refinement.
To that end, we choose the coarsest grid to have $n=20$ points along each direction and refine by a factor of two up to $n=160$ points per direction. We compute ${\rm CFL} = 1/\sqrt{3\mu+\lambda}$ and set the time-step size $k = 0.4 \sqrt{|J|} \cdot {\rm CFL}$, where $J$ is the Jacobian as described in Section \ref{sec::ElWHI_FD}. In Table \ref{tab::ConvergenceFD} we display estimated rates of convergence and observe second order convergence, as expected.

\subsection{Verification of Corrected Time-Steppers}\label{sec::TimeRemoval}
We consider solving the elastic Helmholtz equation with $\lambda = \mu = 1.0$
and choose the forcing so that the exact solution is the same as 
\eqref{eqn::QuarticExample}. We take the frequency $\omega = 1$ and 
enforce homogeneous Dirichlet conditions on the boundary of the 
square $(x,y) \in [0,1]^2$. As the solution
is a fourth degree polynomial, choosing $p = 4$ should ensure that the 
solution to the discrete elastic Helmholtz equation is precisely 
\eqref{eqn::QuarticExample} (up to floating point arithmetic errors) . We use the conjugate gradient accelerated version of El-WaveHoltz with the corrected second order centered time-stepping scheme
presented in Section~\ref{sec::ExplicitTime-Stepping} or with the corrected implicit time-stepping scheme from Section \ref{sec:whimp}.  

We partition the domain into four quadrilaterals
of equal side length $h = 0.5$, set the relative conjugate gradient residual 
tolerance to $10^{-15}$, and consider the error as 
the time-step size is decreased. 

We see from Figure~\ref{fig:time_error_conv} that the standard centered scheme
leads to a discrete solution that converges to the true solution at second order as $\Delta t$ decreases.
The modified schemes on the other hand achieves relative errors that are near machine precision with the error level varying slightly with the time-step size. For the remaining numerical examples unless we state otherwise 
we use the explicit modified time-stepping scheme to remove time discretization errors.

\subsection{Accuracy of the Symmetric Interior Penalty Discontinuous Galerkin Method}

\begin{figure}[ht]
\graphicspath{{figures/}}
\begin{center}
\includegraphics[width=0.49\textwidth,trim={0.0cm -0.1cm 0.0cm 0.0cm},clip]{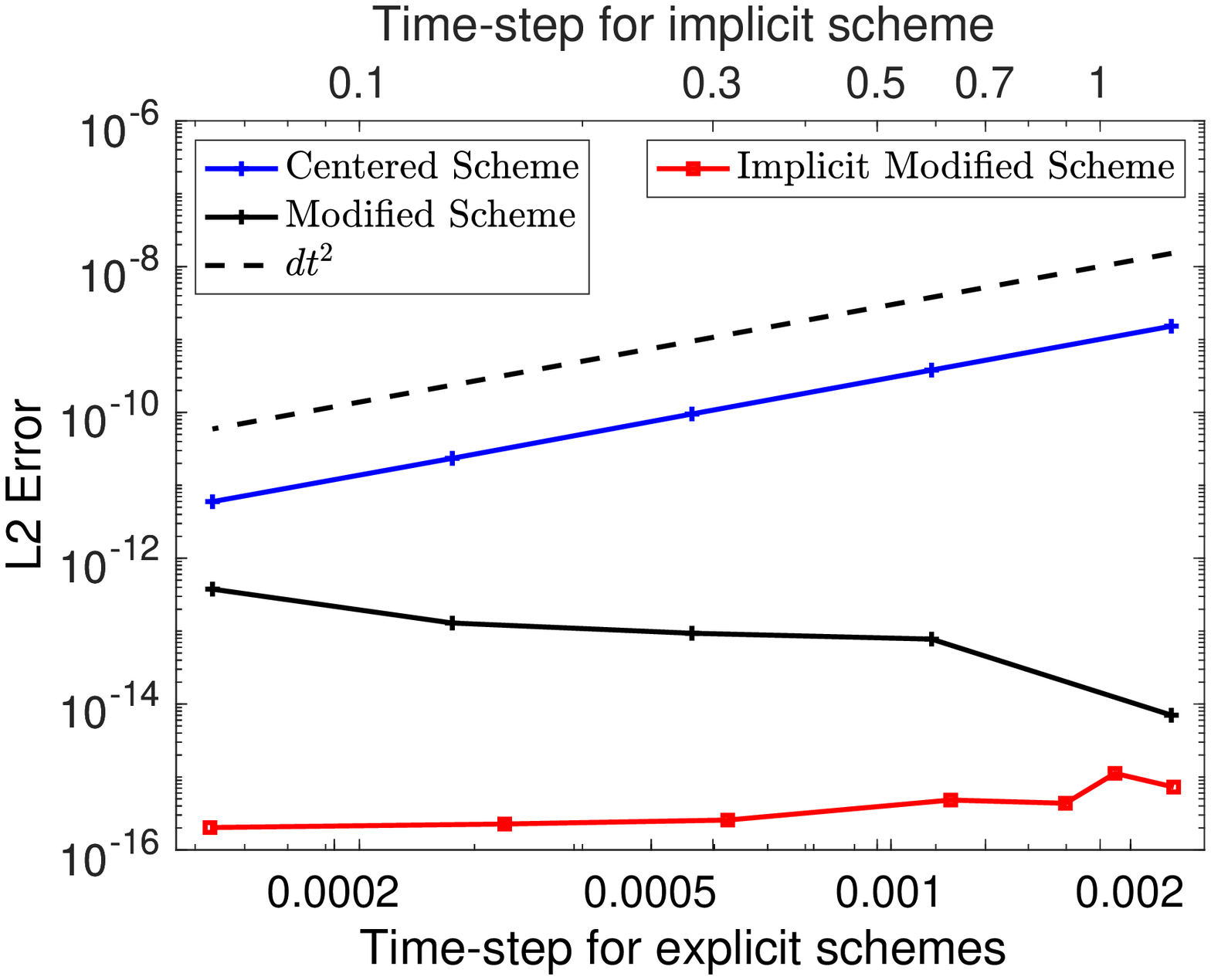} \ \
\includegraphics[width=0.465\textwidth]{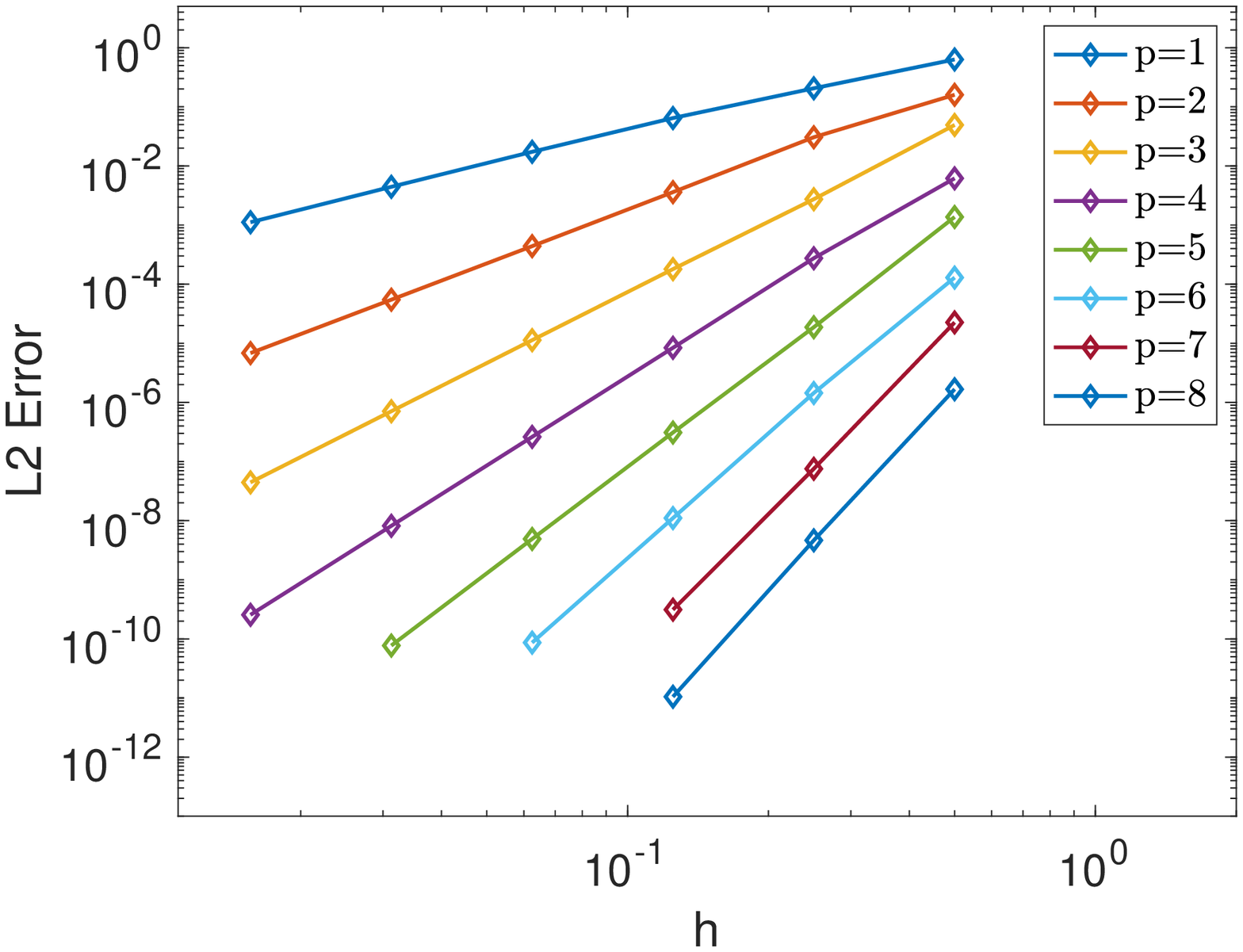}
\caption{(Left) Convergence of the discrete WaveHoltz solution to the true solution of the discrete Helmholtz problem. Note that the time-step size for the explicit and implicit methods are on the bottom and top, respectively, of the figure.  (Right) Convergence of the discrete WaveHoltz solution to the true solution of the discrete Helmholtz problem
for a manufactured solution.  \label{fig:time_error_conv}}
\end{center}
\end{figure}

Next we verify the rates of convergence for our symmetric interior penalty DG solver and for non-homogeneous problems
using an example taken from \cite{el_dg_dath}.
We consider the unit square $S = [0,1]^2$ and impose Dirichlet conditions 
on the boundary. The boundary conditions and forcing are chosen so that 
the Helmholtz solution is
  \begin{align*}
    u(x,y) &= \sin(k_xx + x_0)\sin(k_yy + y_0), \\
    v(x,y) &= -\sin(k_xx + x_0)\sin(k_yy + y_0),
  \end{align*}
where $k_x = 2.5\pi$, $k_x = 2\pi$, $x_0 = 5$, and $y_0 = -10$.
The mesh used is a uniform discretization of the unit square 
split into smaller squares of side-length $h = 1/2^n$
for $n = 1,\dots,6$.

\begin{table}[htb]
\caption{Estimated rates of convergence for the spatial discretization.\label{table::DGOrder}} 
\begin{center} 
 \begin{tabular}{|c|c|c|c|c|c|c|c|c|} 
\hline 
 p & 1 & 2 & 3 & 4 & 5 & 6 & 7 & 8 \\ 
 \hline 
  & 1.84 & 2.94 & 4.00 & 4.94 & 6.01 & 6.86 & 8.07 & 8.64 \\
\hline 
\end{tabular} 
\end{center} 
\end{table} 

We set $\omega = 1$ and choose the material parameters to be the constants $\mu = 1$, $\lambda = 2$. Here we use the modified time-stepping scheme of Section~\ref{sec::ExplicitTime-Stepping} with ${\rm CFL} = 0.4$.

The errors are plotted in Figure~\ref{fig:time_error_conv} as a function of the grid size $h$. We additionally display estimated rates of convergence
calculated using linear least squares 
in Table~\ref{table::DGOrder} from which it is clear that 
the WaveHoltz method converges with optimal 
rates with the error corrected time-stepper (which is formally only second order accurate in $\Delta t$).

\subsection{Effects on Number of Iterations from Number of Periods and Accuracy}
\begin{figure}[ht]
\graphicspath{{figures/}}
\begin{center}
\includegraphics[width=0.95\textwidth,trim={8.4cm 1.0cm 8.4cm 1.0cm},clip]{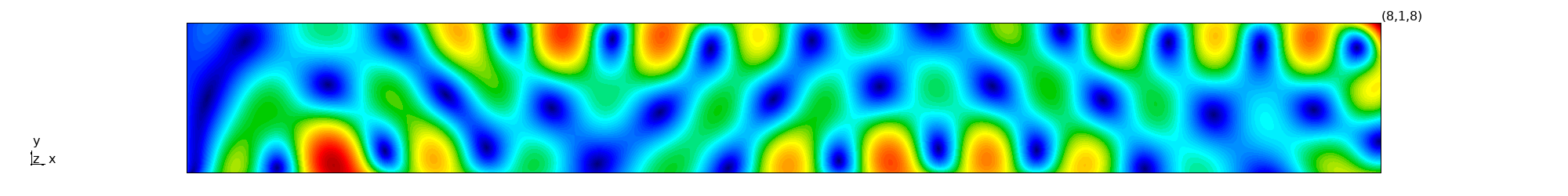}
\includegraphics[width=0.95\textwidth,trim={8.4cm 1.0cm 8.4cm 1.0cm},clip]{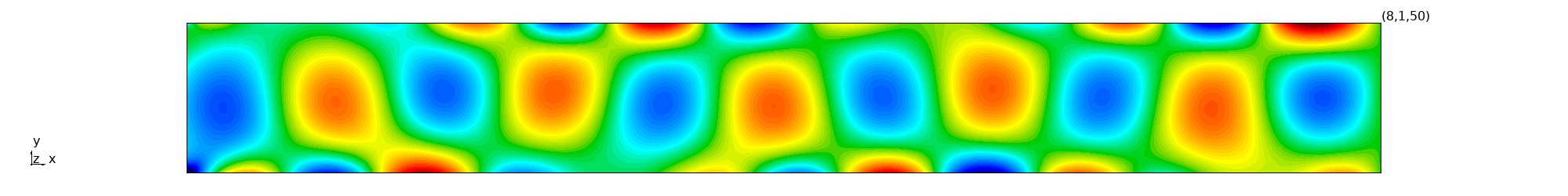}
\includegraphics[width=0.95\textwidth,trim={8.4cm 1.0cm 8.4cm 1.0cm},clip]{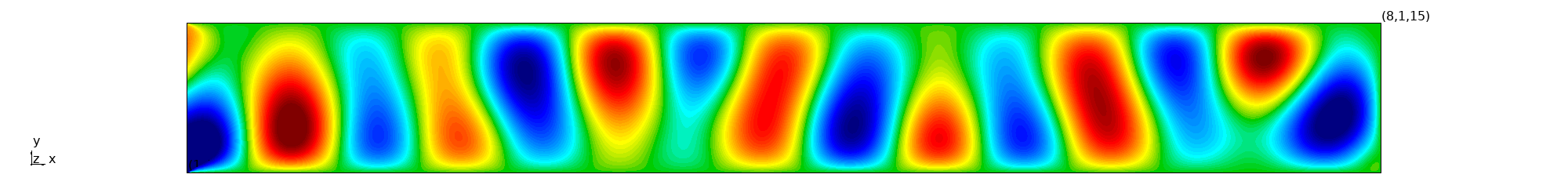}
\includegraphics[width=0.95\textwidth,trim={8.4cm 1.0cm 8.4cm 1.0cm},clip]{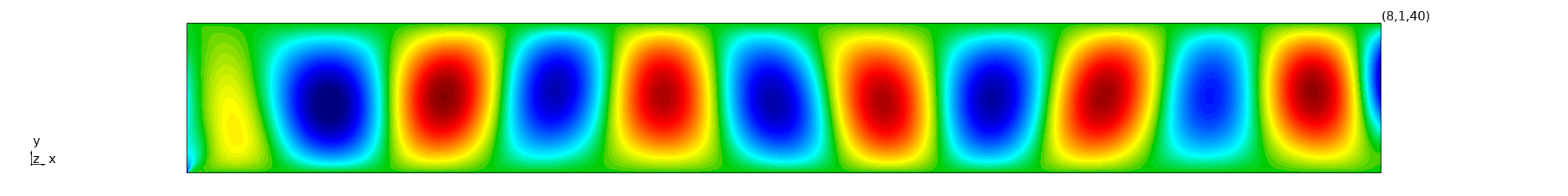}
\caption{From top to bottom: displacement magnitude, $\sigma_{xx}$, $\sigma_{xy}$ and $\sigma_{yy}$. The domain is $[0,8] \times [0, 1]$ and the color scales are $[0,8]$, $[-50,50]$, $[-15,15]$ and $[-40,40]$ respectively. \label{fig:period}}
\end{center}
\end{figure}
In this section we investigate the efficiency of the filter (\ref{eq:many_periods_filter}), defined over $K$ periods, for various values of $K$. Let $N_T$ be the number of time-steps for one period. Then we expect the reduction in the number of all-to-all communications to be $KN_T$ when compared to a direct discretization of 
(\ref{eq:elastic-frequency}).
%(\ref{sec:elastic-wave-eqation-frequency}). 
Here we consider energy conserving boundary conditions, for which we can use the conjugate gradient method and avoid the need to store a Krylov subspace. We note that for problems with impedance or non-reflecting boundary conditions (or with lower order damping terms), El-WaveHoltz will still result in a positive definite but non-symmetric system which can be solved e.g. with GMRES. In that case, we also expect that the size of the GMRES Krylov subspace will decrease by a factor $KN_T$ compared to direct discretization, and by a factor $K$ compared to using a single period for the filter procedure. 

\begin{table}[ht]
 \caption{The table displays the number of iterations required and the efficiency of the longer times to reduce the relative residual by a factor $10^{-10}$ for the two cases (described in the text).\label{tab:period}} 
\begin{center} 
 \begin{tabular}{|c|c|c|c|c|c|c|} 
\hline 
Periods  & 1 & 2 & 3 & 4 & 5 & 10\\
\hline
Case 1 (\#iter) & 124 & 69 &  51 & 43 &  39 & 28\\
Efficiency & 1 &   0.90 &   0.81 &  0.72 &  0.64  &  0.44\\
\hline
Case 2 (\#iter) & 151 & 96 &  78 & 68 &  62 & 53\\
Efficiency & 1 & 0.79 & 0.65 & 0.56 &  0.49 &   0.29\\
\hline 
\end{tabular} 
\end{center} 
\end{table} 

With these obvious advantages of filtering over $K$ periods, it is natural to ask how the number of iterations are affected by the increased filter time. In this experiment we numerically investigate this. To do so, we use the corrected explicit version of the DG solver and consider the shaking of a bar of (unitless) length 8 and height 1. We impose free surface boundary conditions on the top, bottom, and right of the domain, and on the left we set the boundary conditions to be 
\[
u(0,y,t) = v(0,y,t) = \cos(\omega t).
\] 
The base computational mesh uses 8 square elements each with side length 1, which we uniformly refine by dividing each element into 4 parts for some number of refinements. We set $\lambda = 2$, $\mu=1$, $\omega = 5.123$ and ${\rm CFL} = 0.8$. We consider 2 cases: Case 1 uses $p = 5$ and refines the base grid 3 times, Case 2 uses $p=15$ and refines one time. For both cases we use conjugate gradient and count the number of iterations it takes to reduce the relative residual by a factor $10^{-10}$. The solution, along with the components of the stress tensor $\sigma_{xx}$, $\sigma_{xy}$ and $\sigma_{yy}$, are displayed in Figure \ref{fig:period}. The number of iterations for the two cases and the relative efficiency are tabulated in Table \ref{tab:period}. Here the relative efficiency is computed via 
$N_{1}/KN_{K}$ where $N_{j}$ is the number of iterations
required to reach convergence using $j$ periods. As can be seen, the efficiency is relatively high when the number of periods are small and can thus be deployed if the all-to-all communications (or the size of a GMRES Krylov space) becomes a limiting factor.     

\subsection{Iteration Count as a Function of Frequency for Rectangles and Annular Sectors}
For energy conserving boundary conditions, the theoretical prediction (which is also observed experimentally) is that the number of iterations scales as $\omega^{d}$ in $d$-dimensions. In this and the next section, we study how the number of iterations depends on the frequency in two and three dimensions. In this section we additionally investigate the dependence of the number of iterations on frequency for different geometries when the conjugate residual method is used. Here, we study these properties via three different computational domains: a rectangle, a quarter annulus, and a half annulus (all with a characteristic length of 5). We use the finite difference method together with the standard explicit second order time-stepping scheme. For each of the geometries we consider the set of frequencies $\omega = k + \sqrt{2}/10,$ with $ k = 3,4,\dots, 40$. 

\begin{figure}[ht]
\graphicspath{{figures/}}
\begin{center}
\includegraphics[width=0.25\textwidth]{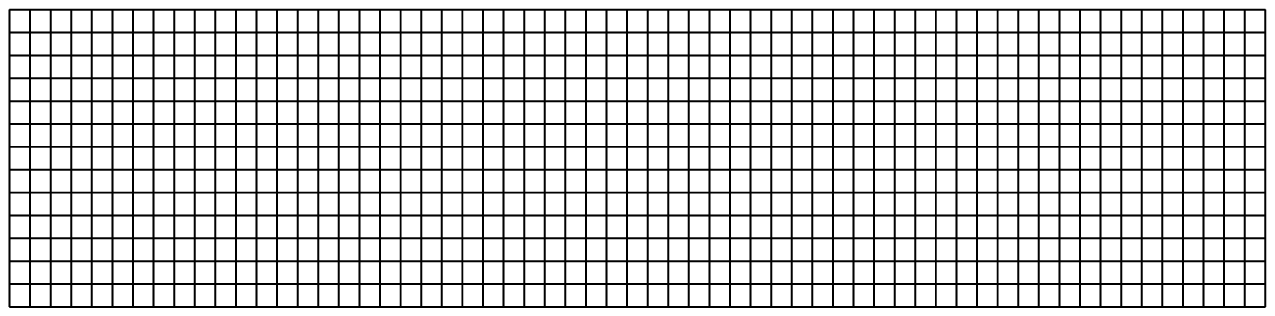}
\includegraphics[width=0.25\textwidth]{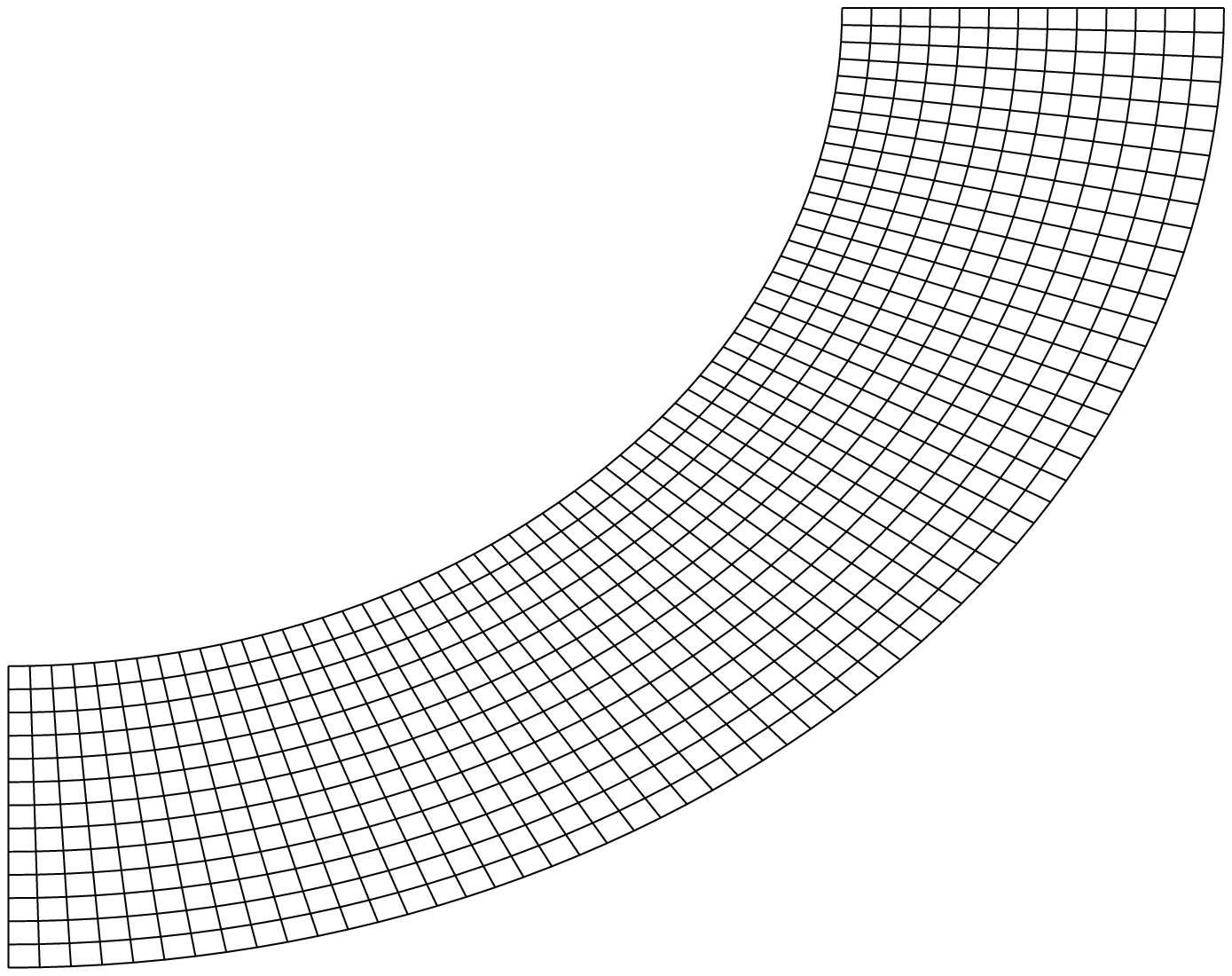}
\includegraphics[width=0.25\textwidth]{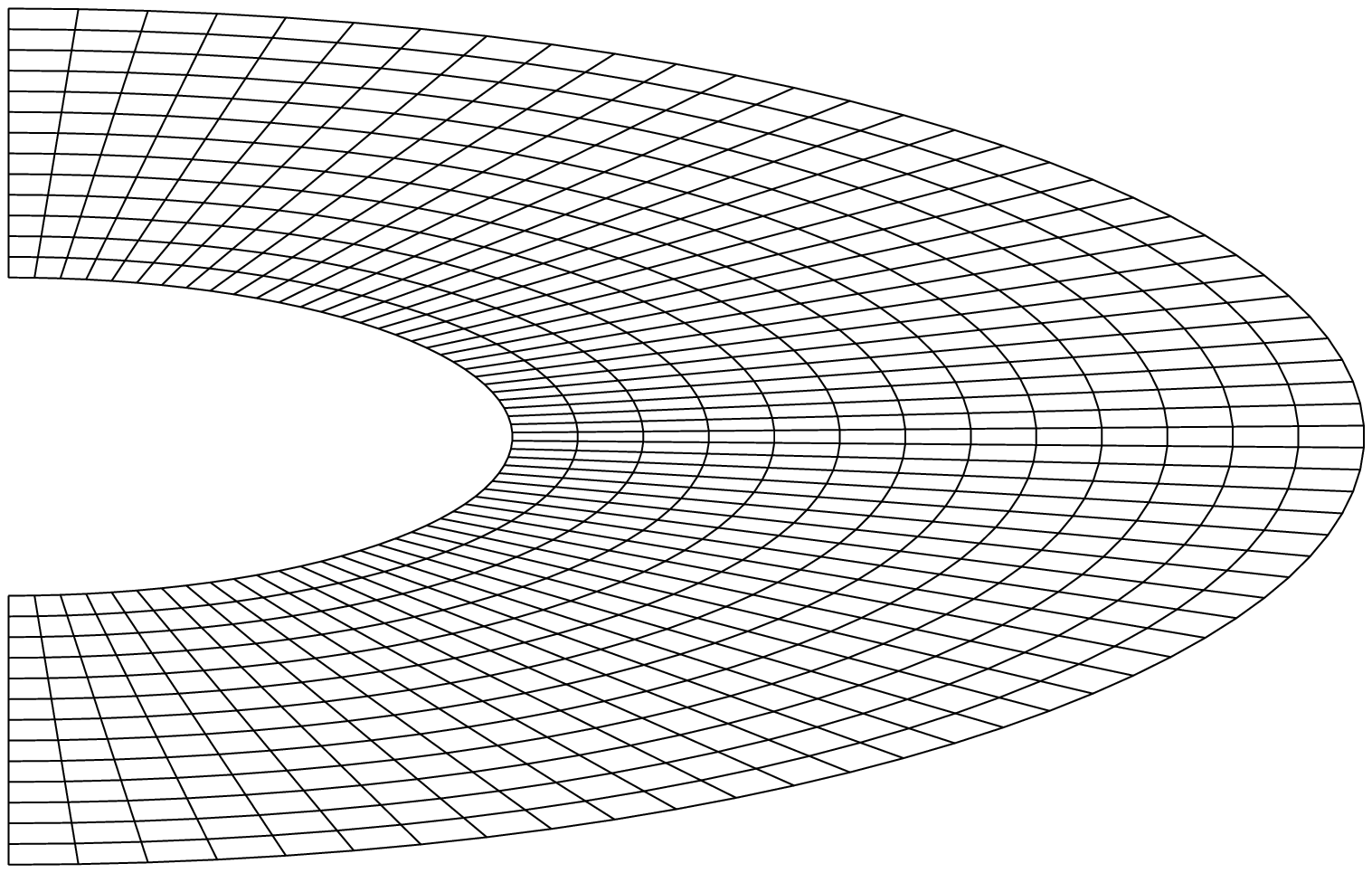}
\caption{A plot of the (Left) rectangle, (Middle) quarter annulus, and (Right) half annulus regions. \label{fig:rectAnnularMesh}}
\end{center}
\end{figure}

Let $q$ and $r$ be the coordinates in the (reference) unit square. We set $n_q$ and $n_r$ to be the number of cells in each coordinate direction. The (spatial) step size is given by $h_q = 1/n_q$ and $h_r = 1/n_r$, and our grid on the unit square is given by 
\begin{align*}
q_i = ih_q ,\ i = 0,\dots, n_q, \ \ r_j = jh_r ,\ j = 0,\dots, n_r.\\
\end{align*} 
We set the arc length of the outer arc at radii $r_{\text{out}}$ of the annular sector to be the length, $L = 5$. Thus for the quarter annulus we have $r_{\rm out} = \frac{2L}{\pi}$, and for the half-annulus we have $r_{\rm out} = \frac{L}{\pi}$. For both cases, we take $r_{\rm in} = r_{\rm out} - 1$. Precisely, the coordinates of the two grids used for the quarter annular sector and the half annular sector are 
\begin{align*}
x_{ij} &= (r_{\rm {in}} + (r_{\rm {out}} - r_{\rm{in}}) q_i)\cos\left(n_{\rm {an}} \frac{\pi}{2} r_j - \frac{\pi}{2} \right),\\
y_{ij} &= (r_{\rm {in}}+ (r_{\rm {out}} - r_{\rm {in}}) q_i)\sin\left(n_{\rm {an}} \frac{\pi}{2} r_j - \frac{\pi}{2} \right).
\end{align*}
Here $n_{\rm an}$ is either 1 or 2 to indicate the quarter or half annulus, respectively,
and plots of each region are displayed in Figure~\ref{fig:rectAnnularMesh}. 
We set $n_r = 4L\omega + 1, n_q = 4\omega +1$ so that the number of points per shear wavelength 
is about 20.
% Insert picture of domain?
\begin{figure}[ht]
\graphicspath{{figures/}}
\begin{center}
\includegraphics[width=0.47\textwidth,trim={0.0cm 0.0cm 0.0cm 0.0cm},clip]{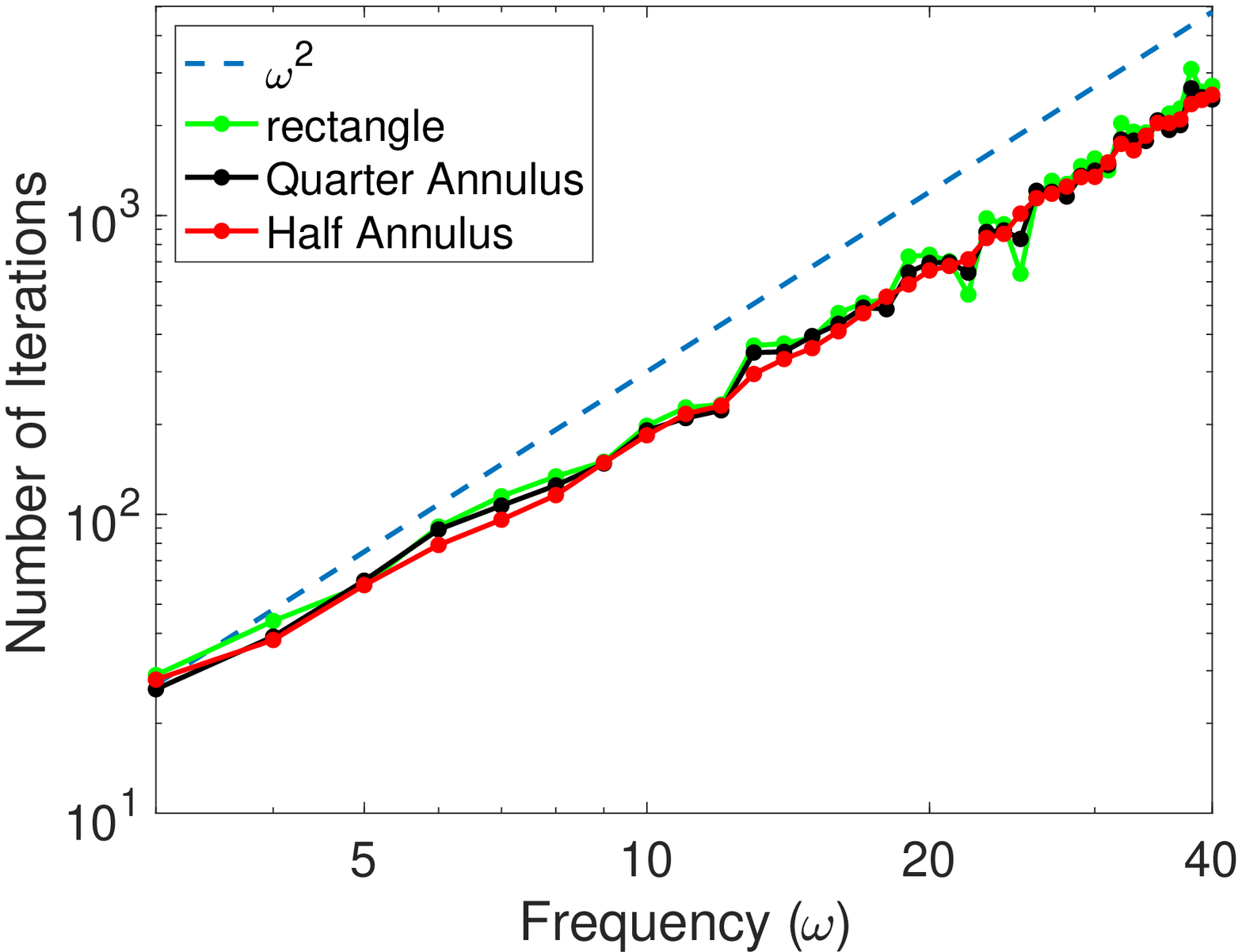}
\includegraphics[width=0.46\textwidth,trim={0.0cm 0.0cm 0.0cm 0.0cm},clip]{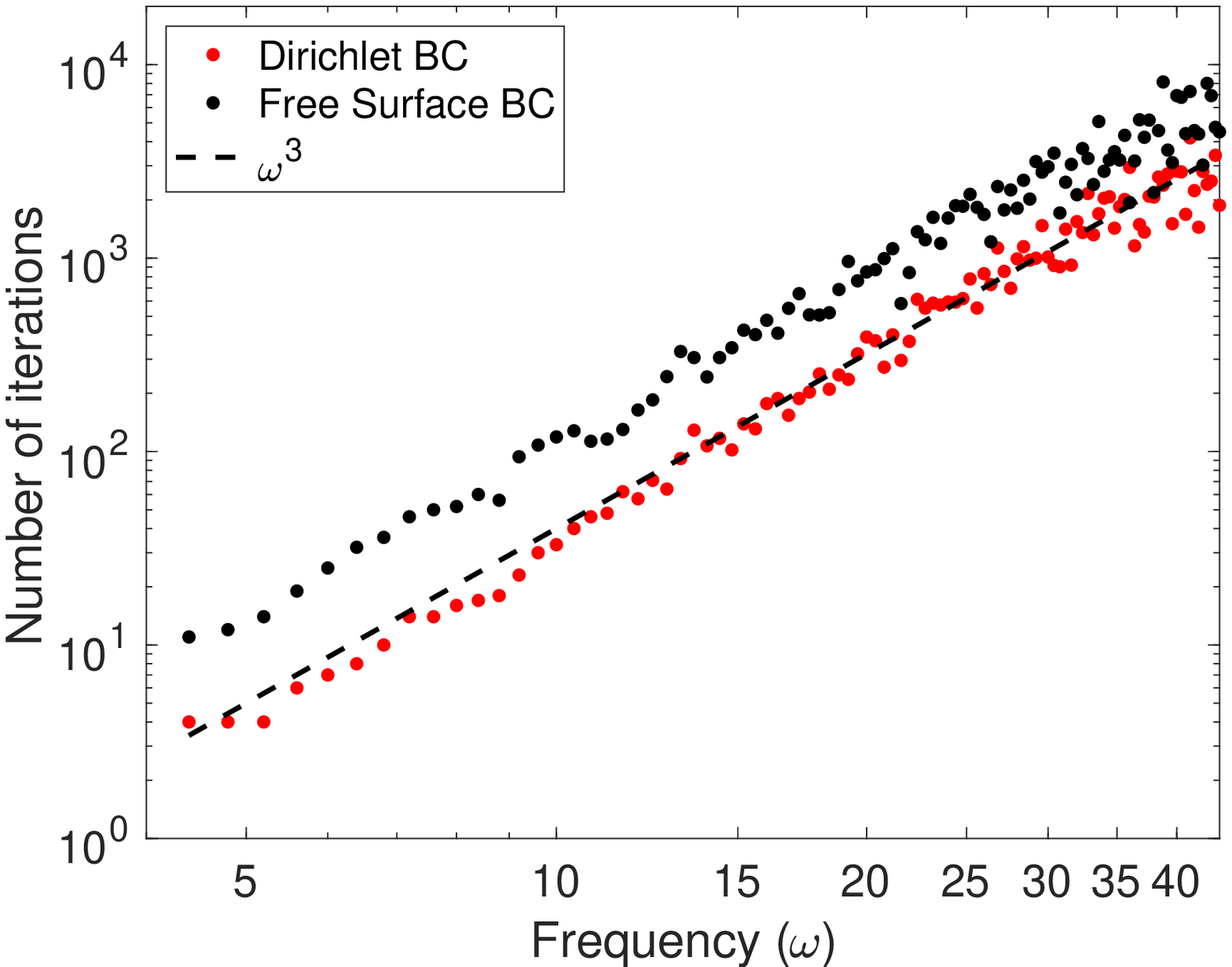}
\caption{The number of iterations as a function of frequency to reach convergence for (Left) a rectangle, quarter circle and half circle, 
and (Right) the unit cube with Dirichlet or free surface conditions. \label{fig:rectAnnular5}}
\end{center}
\end{figure}  

For the forcing we use a discrete approximation of the delta function with amplitude $\omega^2\cos(\omega t)$. We locate this point source at $(x_{i^*j^*},y_{i^*j^*})$ where $i^* = j^* = (n_q+1)/2$ so that it is close to $(0.5,0.5)$ in physical space. In Figure~\ref{fig:rectAnnular5} we display the number of iterations required to reduce the relative residual in the conjugate residual method by a factor $10^{-8}$. From Figure \ref{fig:rectAnnular5} we see that the results for the three geometries are very similar, indicating that (in this example at least) the geometry has little to no effect on the number of iterations needed. Moreover the number of iterations grow as $\omega^2$, as expected. 

\subsection{Effect of Boundary Conditions in a Cube}
In this experiment we consider the unit cube with either Dirichlet boundary conditions on all sides, or with free surface boundary conditions on the top and bottom ($z=0$ and $z = 1$) with Dirichlet boundary conditions on all other sides. We use the 3D version of the finite difference method described above with the standard explicit time-stepping method. Here we use the forcing 
\[
f_j(\bx,t) = A_j e^{-\frac{\sigma}{2} \|\bx - \bx_j\|^2},
\]
with $A_j \sim \sqrt{\sigma^d}$ and with $\sigma \sim \omega$ so that each of the components of the forcing approaches a delta function as $\omega$ grows. We select $\bx_j$ sightly different for each $j$ so that both $\nabla \times \bff \neq 0$  and $\nabla \cdot \bff \neq 0$, resulting in a solution with both shear and pressure waves.  

We use the conjugate residual method, keep the product $h \omega = 0.4 $ fixed, and report the number of iterations required to reduce the initial residual (starting from zero initial data) by a factor $10^{-9}$. The result, which can be found in Figure \ref{fig:rectAnnular5}, confirms the prediction from \cite{WaveHoltz} that the number of iterations scale as $\omega^3$.

\subsection{Iteration Count as a Function of Wave Speed Ratio}
The length of a domain, when measured in number of wavelengths, will increase either if the physical domain size is increased or if the wave speed is reduced. The pressure and shear wave speeds are $C_{\rm p} = \sqrt{(2\mu+\lambda)/\rho}$ and $C_{\rm s} = \sqrt{\mu/\rho}$, respectively.  We expect that the number of iterations will depend on the smallest wave speed, but for El-WaveHoltz there is no intuitive reason to think that a problem with $C_{\rm p} \gg C_{\rm s}$ should be more difficult than a problem with $C_{\rm p} \approx C_{\rm s}$. We note, however, that such behavior has been reported in the literature (see e.g. Table 3.1 on page 11 of \cite{treister2018shifted}) for other methods. 
\begin{table}[ht]
 \caption{The effect on iteration count depending on different combinations of $\lambda$ and $\mu$. \label{tab:poisson}} 
\begin{center} 
 \begin{tabular}{| l |c|c|c|c|c|c|c|c|c|} 
\hline 
$\lambda$ & 1 & 2& 4& 8&  16 & 32 & 64 & 1 & 1\\  
\hline
$\mu$ & 1 & 1& 1& 1&  1 & 1 & 1 & 1/4 & 1/16  \\
\hline
$h_{\rm max} \times 10^{2}$ & 3.37 & 3.37 & 3.37 & 3.37 & 3.37& 3.37& 3.37 & 1.10 & 0.852\\
\hline
$h_{\rm min} \times 10^{2} $& 2.25 & 2.25 & 2.25  & 2.25 & 2.25 & 2.25 & 2.25 & 1.70& 0.547 \\
\hline
\#Iter. & 45 & 47 & 35 & 38 &38 & 56 & 44 & 126 & 247\\
\hline
\#Iter. $\times \sqrt{\mu} $ & 45 & 47 & 35 & 38 &38 & 56 & 44 & 63 & 62\\ 
\hline 
\end{tabular} 
\end{center} 
\end{table} 

\begin{figure}[ht]
\graphicspath{{figures/}}
\begin{center}
\includegraphics[width=0.49\textwidth,trim={6.5cm 3.0cm 0.9cm 2.0cm},clip]{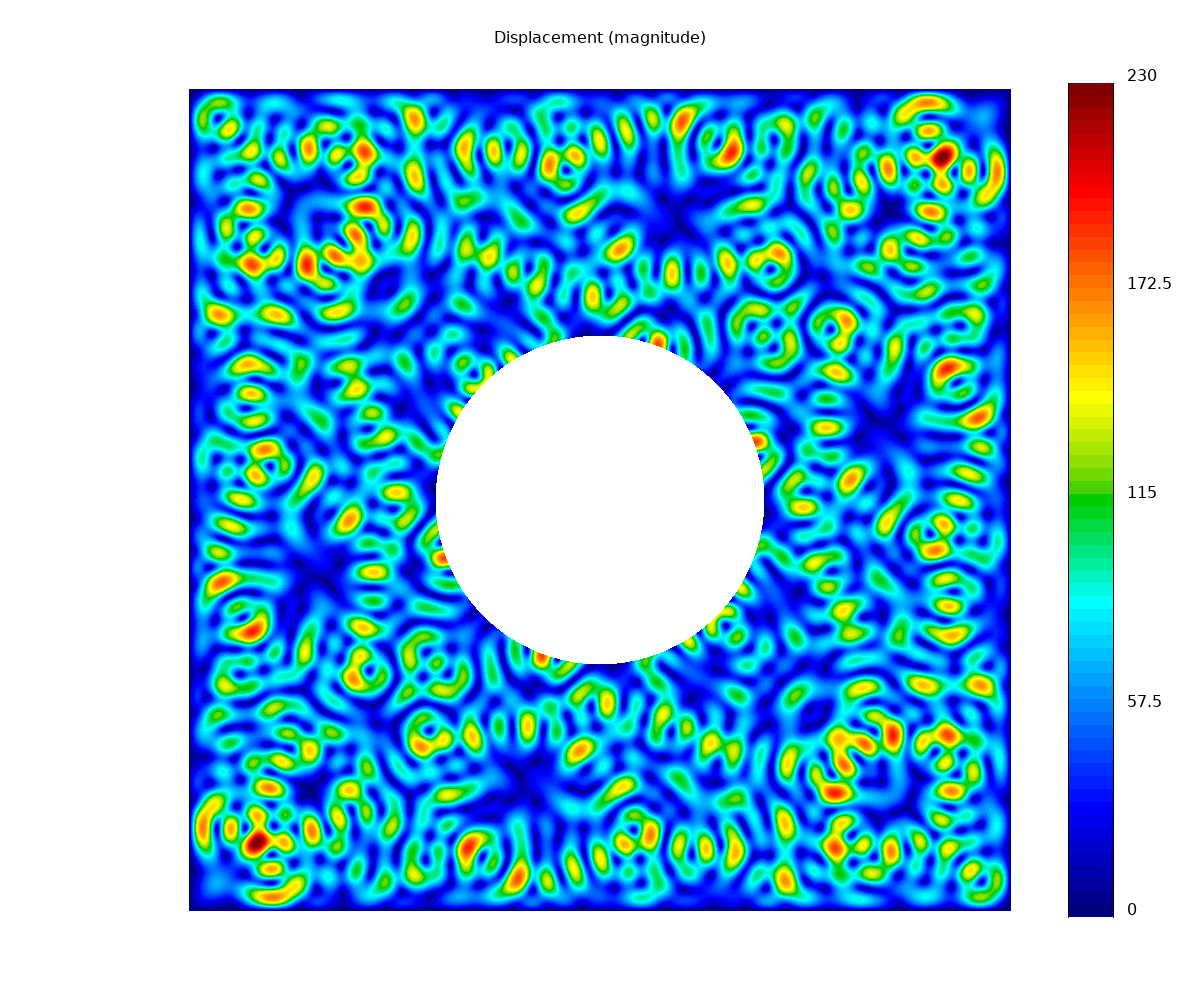}
\includegraphics[width=0.49\textwidth,trim={6.5cm 3.0cm 0.9cm 2.0cm},clip]{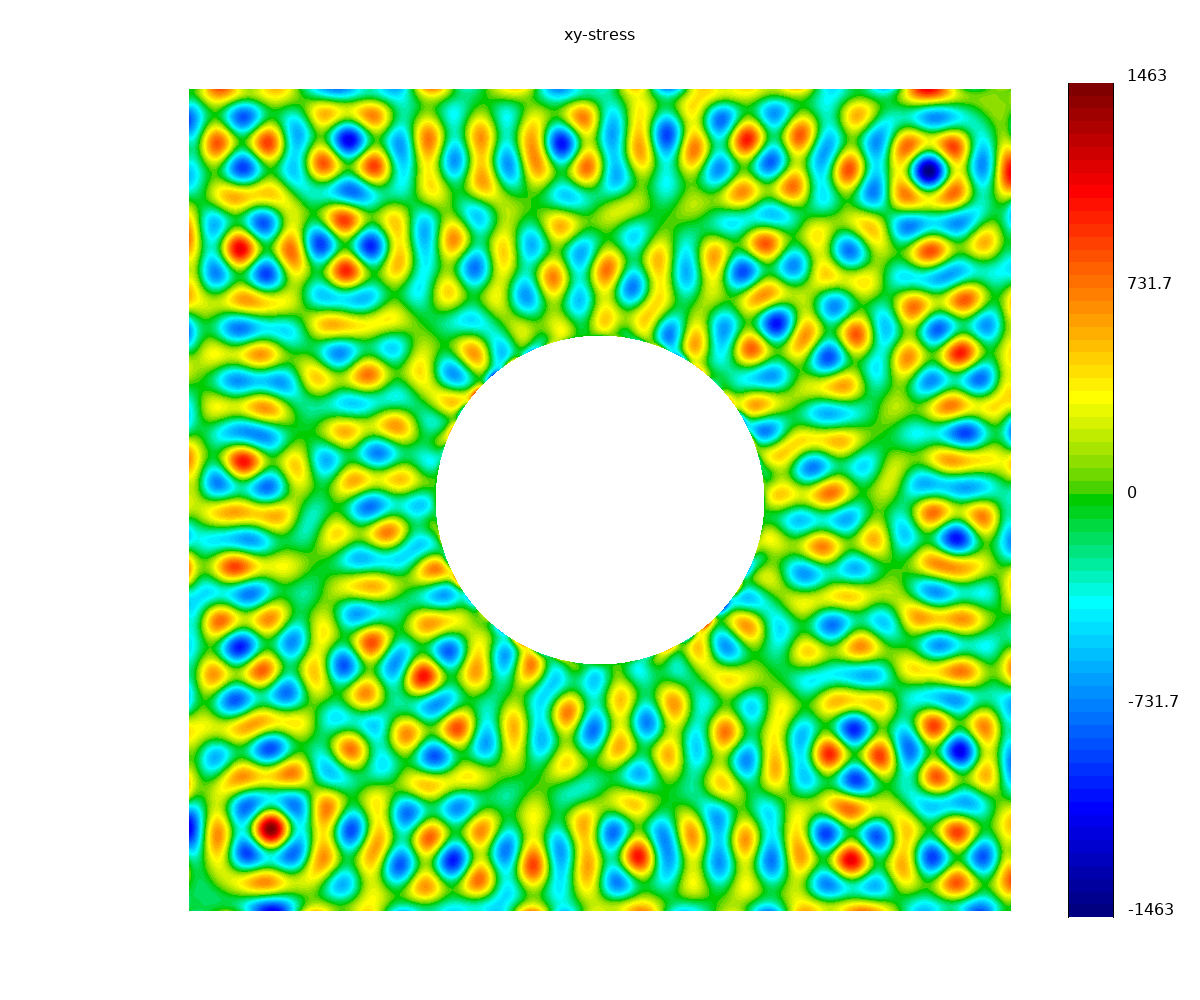}
\includegraphics[width=0.49\textwidth,trim={6.5cm 3.0cm 0.9cm 2.0cm},clip]{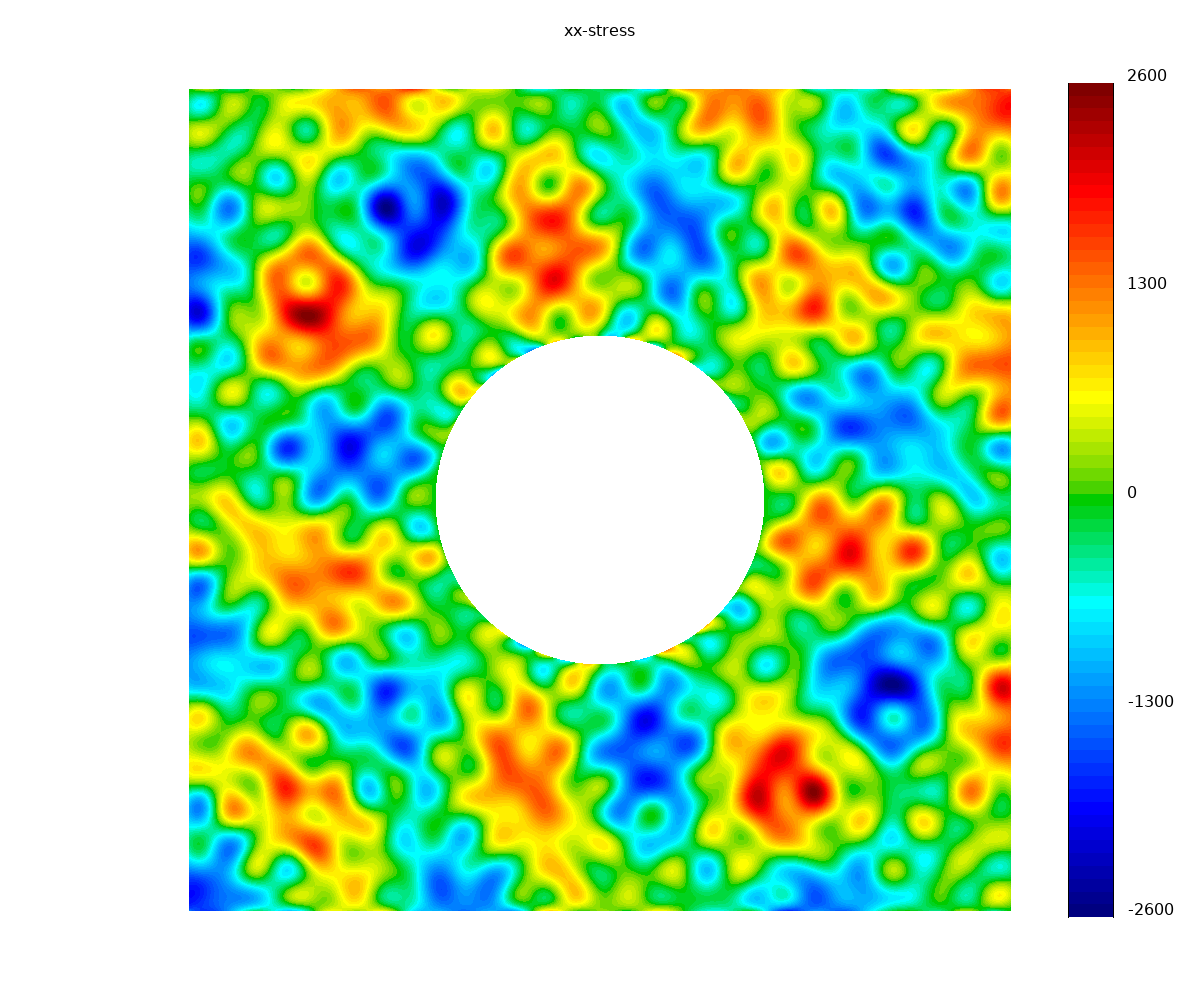}
\includegraphics[width=0.49\textwidth,trim={6.5cm 3.0cm 0.9cm 2.0cm},clip]{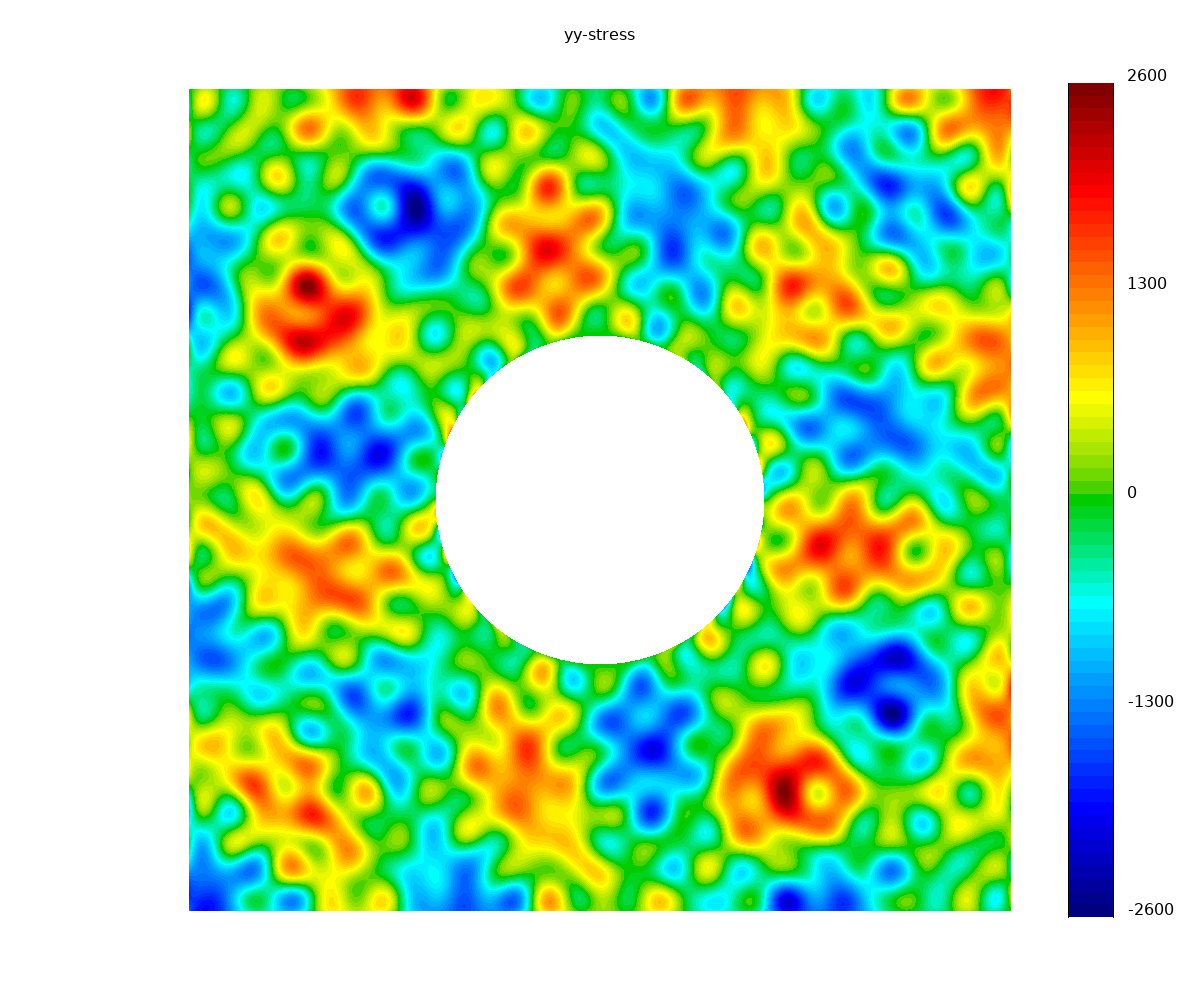}
\caption{From top left to bottom right: displacement magnitude, $\sigma_{xy}$, $\sigma_{xx}$ and $\sigma_{yy}$. \label{fig:poisson}}
\end{center}
\end{figure}

To experimentally investigate how well El-WaveHoltz works for different combinations of $\mu$ and $\lambda$, we use the SIPDG solver with the corrected explicit time-stepper for a geometry consisting of the unit square with a circular hole cut out (see Figure \ref{fig:poisson}). This is the mesh {\tt square-disc-nurbs.mesh}, which is part of the MFEM distribution. The Lam\'e parameters are constant in space and we choose the number of refinements so that the solution is well resolved (the largest and smallest element size is reported in Table \ref{tab:poisson}). We impose the boundary conditions 
\[
u(0,y,t) = v(0,y,t) = \cos(\omega t),
\] 
on the outer part of the domain, and let the circular hole be free of traction. For all experiments we set $\omega = 25.12$, ${\rm CFL} = 0.8$ and we evolve the El-WaveHoltz iteration over $K = 3$ periods.  We stop the CG iteration when the relative residual falls below $10^{-6}$. In Figure \ref{fig:poisson} we display the magnitude of the displacement and the components of the stress tensor $\sigma_{xy}$, $\sigma_{xx}$ and $\sigma_{yy}$ for the case when $\lambda = 1$ and $\mu = 1/16$.   

The results, displayed in Table \ref{tab:poisson}, show that El-WaveHoltz appears to be robust with respect to the ratio between $\lambda$ and $\mu$. Moreover, the number of iterations to reach the desired tolerance is primarily a function of $\mu$, or equivalently, the shear wave speed $C_{\rm s}$.  

\begin{table}[ht]
\caption{Comparison of the number of iterations for the three different methods. \label{tab:HH_vs_WHI_iter}} 
\begin{center} 
 \begin{tabular}{| l | l | r |r|r|r|r|r|r|r|r|r|} 
\hline
&  & $p = 1$& 2 & 3 & 4 & 5 & 6& 7& 8& 9\\
\hline
$h$ &WH  &   57 &   102 &   115 &   126 &   132 &   138 &   141 &   144 &   148 \\
$h/2$ &WH  &   68 &   108 &   114 &   126 &   130 &   133 &   138 &   143 &   146 \\
$h/4$ &WH  &   79 &   107 &   114 &   122 &   130 &   135 &   137 &   142 &   146 \\
$h/8$ &WH  &   83 &   108 &   114 &   125 &   130 &   133 &   137 &   142 &   146 \\
\hline
$h$ &IWH  &   81 &   182 &   160 &   173 &   213 &   237 &   192 &   235 &   201 \\
$h/2$ &IWH  &   97 &   151 &   157 &   168 &   220 &   263 &   188 &   276 &   196 \\
$h/4$ &IWH  &  112 &   147 &   153 &   163 &   171 &   178 &   264 &   188 &   275 \\
$h/8$ &IWH & 117 &   147 &   154 &   165 &   171 &   213 &   217 &     N/A &    N/A \\
\hline
$h$ &HH     & 480 & 25084 & 18298 & 37926 & 80262 & 144863 & 204694 & 230688 & 500000 \\
$h/2$ &HH     &15686 & 32063 & 57338 & 106688 & 256801 & 347279 & 500000 & 500000 & 500000 \\
$h/4$ &HH     &46667 & 79865 & 184331 & 334665 & 500000 & 500000 & 500000 & 500000 & 500000 \\
$h/8$ &HH     &500000 & 304561 & 500000 & 500000 & 500000 & 500000 & 500000 &     N/A &     N/A \\
\hline 
\end{tabular} 
\end{center} 
\end{table}

\subsection{Comparison of Explicit and Implicit El-WaveHoltz with Direct Discretization of the Navier Equations}
In this example we compare the explicit error corrected SIPDG method, the implicit error corrected SIPDG method and the SIPDG method of MFEM's example code \verb+example17p+ extended to the Navier equations (\ref{eq:elastic-frequency}). In this section we will refer to these solvers by the abbreviations WH, IWH and HH respectively. For the Helmholtz SIPDG solver we use GMRES preconditioned by the AMG solver provided by the \verb+HypreBoomerAMG+ class in MFEM. Specifically, we use the MFEM provided implementation of an elasticity AMG solver that incorporates near null-space rigid body modes in the range of interpolation (see \cite{BakKolYanElasticAMG}). The GMRES solver is restarted every 100 iterations. Here we would like to stress that we use the default AMG parameter selections for the HH and IWH solvers and that there could be other parameter selections that would have worked better for the two solvers. Therefore the results below represent a comparison of the ``plain vanilla'' versions of the three methods rather than a comparison between the best tuned methods. We feel that this is gives a baseline comparison of the methods and also represent a typical practical situation where the solver is only one part of a larger problem (say an inverse problem or a design problem) and there is only limited time to tune the solver.

For the implicit solver we must also invert the elasticity operator (but with a shift that preserves its positive definiteness) and we do this using conjugate gradient preconditioned by the same AMG setup as for the Helmholtz SIPDG solver. We always use 10 time-steps for the implicit solver, and for the explicit solver we use ${\rm CFL} = 1.1$ in \eqref{eqn::DG-CFL}. For the WH and IWH solver we solve the El-WaveHoltz problem with conjugate gradient. For all three solvers the tolerance is set to be $10^{-10}$.The implicit solver computations were carried out on a dual Intel Xeon Platinum 8268 2.98 GHz 24-core processor.

We solve the equations on the unit square with a smooth (but narrow) forcing 
\[
\bff = -\frac{200 \omega^2}{\pi} e^{-1.2\omega^2[(x-0.25)^2+(y-0.25)^2]} \left( \begin{array}{c}
-y+0.5 \\
x-0.5
\end{array} \right),
\]
and with free surface boundary conditions on all sides. We consider different refinements and polynomial degrees from 2 to 9 and estimate the error in the solution by computing a reference solution using degree 11 polynomials (note that the solutions are the same up to the tolerance of the iterative solvers since we have eliminated the time errors).

For all of the computations we record the number of iterations (a maximum of 500000) and list them in Table \ref{tab:HH_vs_WHI_iter}. We note that entries for which the simulation could not be run have been marked with N/A for their iteration counts. It can be seen that the fewest number of iterations are achieved with the WH method. We also note that the number of iterations for WH is insensitive to the mesh resolution and has a weak dependence on the polynomial degree. Again, the latter is due to the fact that the linear system we are solving comes from a bounded operator so that the condition number does not depend on $h$. HH, which discretizes an unbounded and indefinite operator, behaves radically different with the number of iterations increasing rapidly with decreasing mesh resolution. In addition, the number of iterations for the HH method increases very quickly with the polynomial degree and, as a result, many of the accurate test cases fail to converge. Similar to the WH method, the IWH method has an iteration count that is relatively robust under grid refinement but with a slight increase with order. We note that it does appear the higher order methods have more variation in iteration count than the lower order methods, and in general the iteration count is larger than for the WH method. 

The number of iterations displayed in Table \ref{tab:HH_vs_WHI_iter} are not useful for comparing the different methods as each iteration comes with a different computational cost. In Table \ref{tab:HH_vs_WHI_rhs} we instead list the number of right hand side (rhs) evaluations. By a right hand side evaluation we mean a single application of the matrix corresponding to the matrix discretizing the elastic operator. For the explicit method, the total number of rhs evaluations is $N_T N_\text{iter}$, where $N_T$ is the number of time-steps needed to evolve the elastic wave equation one period and $N_\text{iter}$ is the number of iterations. For the IWH method we always take $N_T = 10$ so that the number of rhs evaluations is $10 N_\text{inner} N_\text{iter}$, where $N_\text{inner}$ is the number of inner iterations used by the AMG preconditioner (as reported in Table \ref{tab:HH_vs_WHI_iter}). For the HH method the number of rhs evaluations is equal to the number of GMRES iterations.

\begin{table}[ht]
\caption{The number of right hand side evaluations (estimated) for the three different methods. The top four rows display the actual number of right hand side evaluations and the rows below indicate how many times more the HH and IWH method evaluates the right hand side. An infinity sign indicates that the computation did not converge.\label{tab:HH_vs_WHI_rhs}} 
\begin{center} 
 \begin{tabular}{| l | l | r |r|r|r|r|r|r|r|r|r|} 
\hline
&  & $p = 1$& 2 & 3 & 4 & 5 & 6& 7& 8& 9\\
\hline
$h$ &WH  &1425  &  4998  &  9315  & 15120  & 22176  & 30774  & 40326  & 51552  & 64676  \\ 
$h/2$ &WH & 3400  & 10476  & 18354  & 30240  & 43550  & 59318  & 78936  & 102245  & 127458  \\ 
$h/4$ &WH &7821  & 20758  & 36594  & 58438  & 86970  & 120285  & 156728  & 202918  & 254916  \\ 
$h/8$ &WH & 16434  & 41904  & 73188  & 119750  & 173940  & 236873  & 313456  & 405836  & 509686  \\ 
\hline
$h$ &IWH   &  81  &    96  &    67  &    63  &    72  &    72  &    53  &    60  &    48  \\ 
$h/2$ &IWH  &   72  &    67  &    59  &    53  &    69  &    70  &    46  &    59  &    40  \\ 
$h/4$ &IWH   &  59  &    59  &    51  &    48  &    46  &    41  &    57  &    36  &    49  \\ 
$h/8$ &IWH    & 50  &    53  &    47  &    43  &    42  &    46  &    42  &     N/A  &     N/A  \\ 
\hline 
$h$ &HH     & 1  &     5  &     2  &     3  &     4  &     5  &     5  &     4  &     $\infty$  \\ 
$h/2$ &HH    &  5  &     3  &     3  &     4  &     6  &     6  &     $\infty$  &     $\infty$  &    $\infty$  \\ 
$h/4$ &HH     & 6  &     4  &     5  &     6  &     $\infty$  &     $\infty$  &     $\infty$  &     $\infty$  &     $\infty$  \\ 
$h/8$ &HH    &  $\infty$  &     7  &     $\infty$  &     $\infty$  &     $\infty$  &     $\infty$  &     $\infty$ &     N/A  &     N/A  \\ 
\hline
\end{tabular} 
\end{center} 
\end{table} 

As can be seen in Table \ref{tab:HH_vs_WHI_rhs}, the WH method is also more efficient with respect to the number of rhs evaluations (note that we report total number of rhs for the WH method and multipliers for the other methods). The advantage of the explicit method over the implicit method appears to be decreasing with increased accuracy -- both in terms of decreasing mesh size and increased polynomial order. This is not unexpected as the number of time-steps needed for the explicit method grows linearly with the reciprocal of the mesh size, and quadratically with the polynomial degree while the implicit method maintains the number of time-steps constant. In terms of rhs evaluations, the gap between the WH method and the HH method is smaller than between WH and IWH; though the HH method degrades with increasing mesh refinement.
\begin{table}[ht]
\caption{The table reports how many times longer a computation with the HH and IWH method takes compared to the explicit WH method. \label{tab:HH_vs_WHI_time}} 
\begin{center} 
 \begin{tabular}{| l | l | r |r|r|r|r|r|r|r|r|r|} 
\hline
& $p$ /meth & 1& 2 & 3 & 4 & 5 & 6& 7& 8& 9\\
\hline
$h$ &IWH &192 & 197 & 137 & 113 & 103 & 90 & 59 & 73 & 44 \\
$h/2$& IWH &169 & 185 & 125 & 112 & 107 & 121 & 55 & 42 & 29 \\
$h/4$& IWH &261 & 295 & 61 & 152 & 62 & 68 & 71 & 36 & 42 \\
$h/8$& IWH &28 & 30 & 100 & 69 & 62 & 46 & 34 &    N/A&     N/A \\ 
\hline 
 \hline 
$h$ &HH &15 & 164 & 56 & 56 & 60 & 64 & 62 & 59 &     $\infty$ \\
$h/2$& HH&265 & 157 & 109 & 96 & 117 & 113 &    $\infty$&    $\infty$&     $\infty$ \\ 
$h/4$& HH&587 & 311 & 91 & 260 &    $\infty$&    $\infty$&    $\infty$&    $\infty$&     $\infty$ \\
$h/8$ &HH&  $\infty$& 224 &    $\infty$&    $\infty$&    $\infty$&    $\infty$&    $\infty$&    N/A&   N/A \\
\hline 

\end{tabular} 
\end{center} 
\end{table} 

Finally, in Table \ref{tab:HH_vs_WHI_time} we report the increase in wall-clock time (as a multiplicative factor) for the IWH and HH methods relative to the time required to solve the same problem with the explicit WH method. Throughout, the WH method is one to two orders of magnitude faster than the other methods. The HH method becomes less attractive (especially when it stops converging) as the accuracy is increased, while the IWH method improves with increased accuracy. It is of note that neither the number of iterations nor the number of rhs is a good predictor of compute time. Possible causes for this discrepancy are, a) that we only counted one right hand side evaluation per iteration and neglected the cost of the AMG preconditioner, and b) that the GMRES solve for the HH method actually has a significantly higher cost than the conjugate gradient method.

\begin{rem}
Given the additional computational cost per time-step of the implicit scheme, 
it would perhaps be natural to hope for a significant improvement in run-time 
compared to the explicit scheme given how few time-steps are used in comparison.
As briefly discussed in Remark~\ref{rem::Implicit_scheme}, a larger time-step size 
can lead to a less accurate discrete filter. Informally, a coarser grid in time yields a 
wider discrete filter which leads to more frequencies appearing to be ``closer'' to 
resonance than a continuous analysis would suggest. For the above examples with 
uniform refinement, this requires a larger number of iterations to reach convergence
with a corresponding increase in rhs evaluation count as well as 
longer run-times. As we will see in the next section, however, this may no longer be the case
for meshes that are highly refined or have a wide variation in element sizes due to the 
stricter CFL condition for explicit schemes. 
\end{rem}

\subsection{Geometric Stiffness - Quarter Circle with a Small Cutout}
\begin{figure}[ht]
\graphicspath{{figures/}}
\begin{center}
\includegraphics[width=0.3\textwidth,trim={0.0cm 0.0cm 0.0cm 0.0cm},clip]{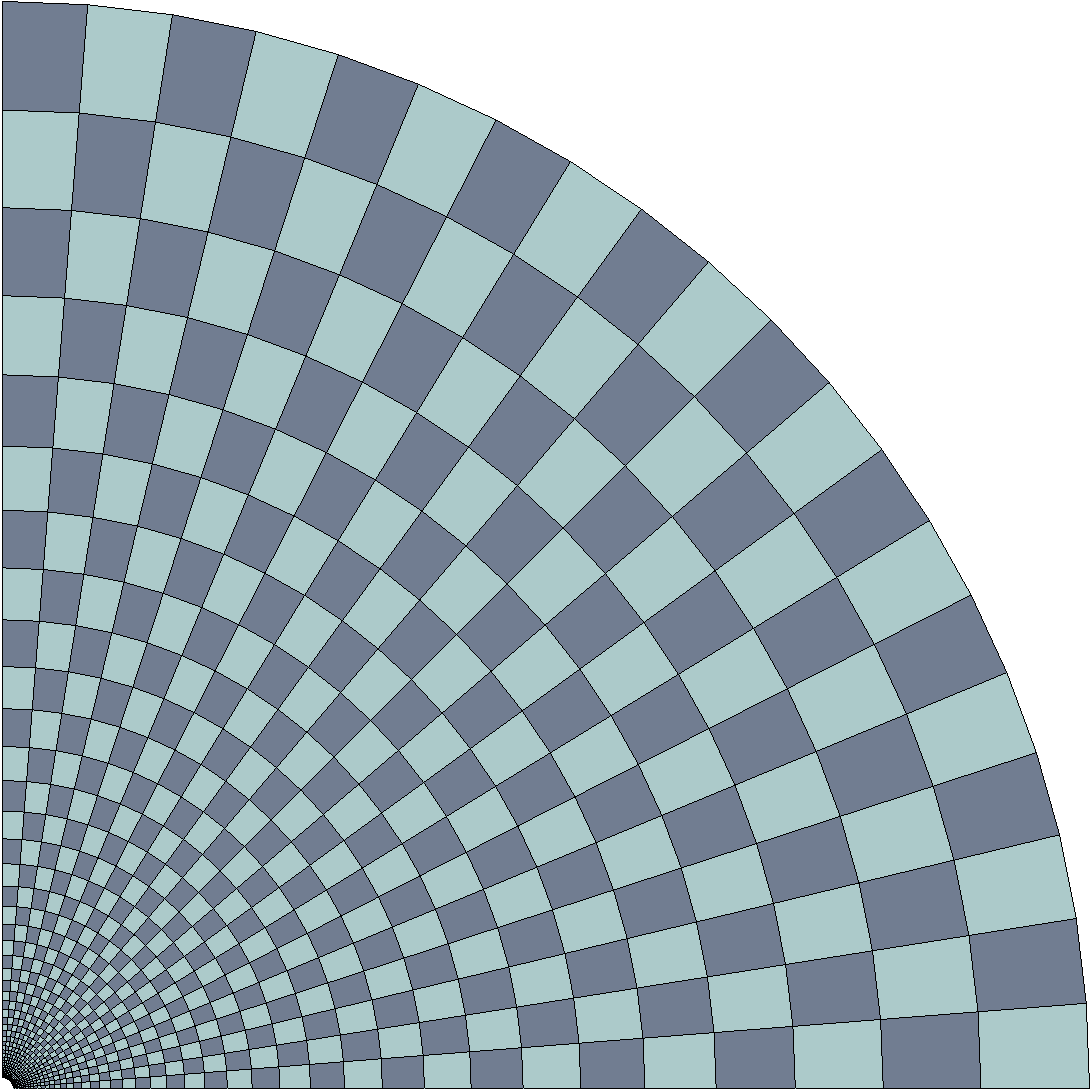}\ \ \ \ \ 
\includegraphics[width=0.3\textwidth,trim={0.0cm 0.0cm 0.0cm 0.0cm},clip]{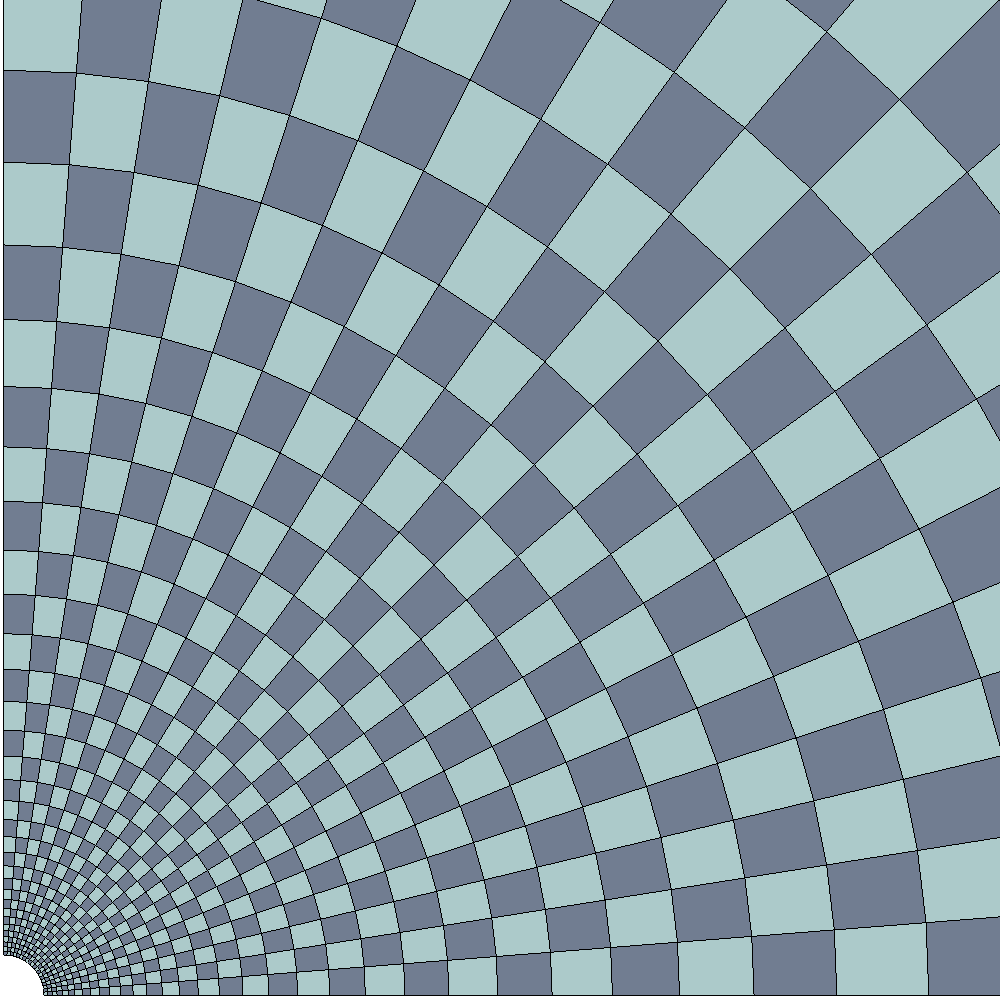}
\caption{Geometry and mesh with a zoom in of the mesh near the semicircle with radii 0.01 to the right. The outer boundary is the unit circle. \label{fig:quart}}
\end{center}
\end{figure}

In this example we highlight the use of the implicit method for meshes with large disparities in element size. In particular we consider the geometry depicted in Figure \ref{fig:quart} which has been discretized with a polar grid with equiangular spacing but with progressively coarser spacings in the radial direction. The ratio between the largest and smallest element side is around one hundred. We impose homogenous Dirichlet boundary conditions on all sides and force the problem with 
\[
\bff = \frac{20000 \omega^2}{\pi} e^{-120\omega^2[(x-0.5)^2+(y-0.5)^2]} \left( \begin{array}{c}
1 \\
1
\end{array} \right).
\]

We again use the SIPDG solver with degree one elements. The implicit solver uses MFEM's PCG solver preconditioned by the AMG solver provided by MFEM's \verb+HypreBoomerAMG+ class. Here we do not use the elasticity specific options (see also \cite{BakKolYanElasticAMG}) as those appeared to make the convergence worse for this problem. For this example we take $\omega = 20$ and $\lambda = \mu = 1$. 
\begin{table}[ht]
\caption{The table compares time spent for the implicit and the explicit method to solve the quarter circle problem.  \label{tab:quart1}} 
\begin{center} 
 \begin{tabular}{| l | r|r|r|r|r|r|r|r|r|} 
%dof  / time
\hline
D.O.F &7200 &  28800 & 115200 & 460800 &1843200 \\
\hline
Time WH &1.00 &   3.5 &  18.3 & 127 &5243 \\
Time WHI-6 (MINR.) & 1.8 &   2.3 &   5.4 &  24.7 & 535 \\
Time WHI-10 &1.2 &   2.6 &   5.3 &  21 & 329 \\
Time WHI-20 & 1.5 &   2.7 &   5.3 &  18.3 & 297 \\
\hline
Time increase WH & -&3.5 & 5.3 & 7.0 &41.1 \\
Time increase WHI-6 (MINR.) &-& 1.3 & 2.3 & 4.6 &21.7 \\
Time increase WHI-10 &-& 2.2 & 2.0 & 3.9 &15.8 \\
Time increase WHI-20 &-&1.9 & 2.0 & 3.5 &16.2 \\
\hline
\end{tabular} 
\end{center}
\end{table} 

For a sequence of 5 uniform refinements we compare the wall-clock time, number of iterations, number of time-steps and (for the three implicit methods) the number of AMG iterations per time-step for the explicit method and the implicit method with 6, 10 and 20 time-steps.  For the 6 time-step solver we do observe that the iteration matrix becomes very slightly indefinite as discussed in the appendix and thus we use MINRES. For the 10 and 20 time-step we don't observe any loss of positive definiteness and use CG to accelerate the iteration. 

The results can be found in Tables \ref{tab:quart1} and \ref{tab:quart2}. In Table \ref{tab:quart1} we display the time (normalized by the time for the explicit method on the coarsest grid) for the different methods along with the number of degrees of freedom for each computation. We also display the time increase factor between two subsequent meshes for each method. It is clear from the table that the implicit method can outperform the explicit method as the mesh is refined. The implicit method that uses MINRES appears to be slightly less efficient.  

With each refinement the number of degrees of freedom increase by a factor of four and the element sides are reduced by a factor of two. This implies that the number of time-steps required for the explicit method also doubles with each refinement. One may therefore expect that the time for the explicit method will grow by a factor of 8 with each refinement. On the other hand, if one makes the assumption that the cost of the multigrid solver scales linearly with the number of degrees of freedom, the cost of the implicit solver should increase by a factor of 4. The results in Table  \ref{tab:quart1} are somewhat consistent with these estimates, with the time increases approaching 8 and 4 from below for the 4 first levels of refinement. The major outlier is the finest mesh for which the wall-clock times increases dramatically. Turning to Table \ref{tab:quart2} we see that the increase in number of AMG iterations per time-step is gradual and does not change drastically between the refinement levels, including the finest level. With this in mind, and as all the methods had a jump in the wall-clock time, we attribute the increase in time to hardware limitations. In the final row of Table \ref{tab:quart2} we compare the the ratio of number of time-steps 
taken between the explicit method and the implicit method with 6, 10 and 20 time-steps. As expected the implicit methods take a significantly fewer number of time-steps relative to explicit methods, especially for highly refined meshes.

\begin{table}[ht]
\caption{The table reports the cost of the IWH method in terms of AMG iterations per time-step for 6, 10 and 20
time-steps. \label{tab:quart2}} 
\begin{center} 
 \begin{tabular}{| l | r|r|r|r|r|r|r|r|r|} 
\hline
D.O.F &7200 &  28800 & 115200 & 460800 &1843200 \\
\hline
AMG / TS WHI-6 (MINR.) & 44.70 & 48.10 & 59.50 & 82.80 &118.30 \\
AMG / TS WHI-10 & 46.11 & 49.09 & 61.50 & 86.00 &129.71 \\
AMG / TS WHI-20 &38.30 & 42.80 & 52.90 & 72.00 &104.76 \\
\hline
AMG increase WHI-6 (MINR.) &-&1.08 & 1.24 & 1.39 & 1.43 \\
AMG increase WHI-10&-& 1.06 & 1.25 & 1.40 & 1.51 \\
AMG increase WHI-20&-& 1.12 & 1.24 & 1.36 & 1.46 \\
\hline
Fewer TS WHI-6 (MINR.) & 1588 & 3212 & 6463 & 3890 &31164 \\
Fewer TS WHI-10 & 953 & 1928 & 3878 & 7779 &15582 \\
Fewer TS WHI-20 &476 &  964 & 1939 & 3890 & 7791 \\
\hline
\end{tabular} 
\end{center} 
\end{table}

\subsection{Materials with Spatially Varying Properties}
We next consider an example taken from \cite{LeVeque:2002kk}
with elastic propagation in a heterogeneous medium.
We define the domain as $\Omega = [0,2]\times[0,1]$
where there is an embedded inclusion, 
$\Omega_I = [0.5,1.5]\times[0.4,0.6]$,
in the middle from a stiffer material. We let 
$\Omega_I$ be a material with $\lambda = 200, \mu = 100$, 
while the domain $\Omega \setminus \Omega_I$ has 
$\lambda = 2$ and $\mu = 1$.
We impose traction-free boundary conditions at $y=0,1$ 
and at $x=1$. We additionally have
  \begin{align*}
    \bu(0,y,t) = \begin{pmatrix}
      0 \\ \cos(\omega t)
    \end{pmatrix},
    \quad 
    \bff(x,y,t)
    \begin{pmatrix}
      0 \\ \delta(|x-x_0| + |y-y_0|)\cos(\omega t)
    \end{pmatrix},
  \end{align*}
where $\delta$ is a delta function centered at $x_0 = 0.1$
and $y_0 = 0.5$.

\begin{figure}[ht]
\graphicspath{{figures/}}
\begin{center}
\includegraphics[width=0.48\textwidth]{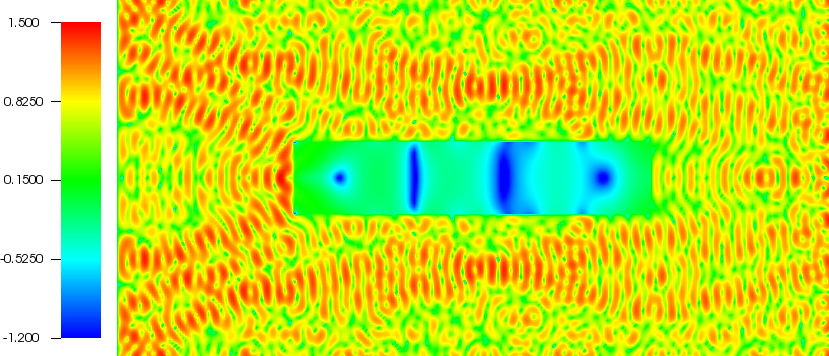}
\includegraphics[width=0.48\textwidth]{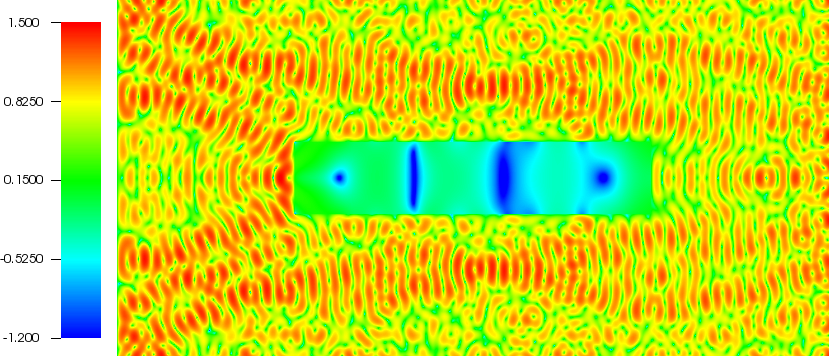}
\caption{The $\log_{10}$ of the magnitude of the displacements for the 
CG accelerated solution of WH for the inclusion problem using sixth order
polynomials within each element. (Left) Solution using a 
grid resolution of at least one element per wavelength, and (Right) two elements per wavelength. 
\label{fig:inclusion_mag}}
\end{center}
\end{figure}

We use the SIPDG method of Section~\ref{sec::SIPDG} with a uniform
quadrilateral mesh and  polynomial degree  $p=6$. We consider two 
meshes with element widths of $h_1 = 1/20$ and $h_2 = 1/40$,
respectively. We use ${\rm CFL} = 0.4$ with the (corrected) explicit 
second order time-stepper of Section~\ref{sec::ExplicitTime-Stepping}, 
integrate over five periods, and accelerate convergence with the
conjugate gradient method with a relative residual tolerance of 
$10^{-5}$. We choose the frequency $\omega = 100$ so that 
we have at least one element per wavelength when using 
element widths of $h_1$, and at least two when the widths are $h_2$. 
We plot the $\text{log}_{10}$ of the magnitude of the 
displacement vector in Figure~\ref{fig:inclusion_mag}.

From Figure~\ref{fig:inclusion_mag} it is clear that the 
solution using one element per wavelength is visually quite 
similar to that of the refined mesh with two elements per wavelength.
Thus using one element per wavelength with a higher order polynomial order 
is sufficient to produce reasonable results, 
as was similarly observed in \cite{tsuji2014sweeping}. 
Moreover, the iteration counts are similar with 634 and 
630 iterations required to reach a relative residual tolerance of 
$10^{-5}$ for elements of size $h_1 = 1/20$ and $h_2 = 1/40$, 
respectively.

\subsection{Vibrations of a Toroidal Shell}
Finally, as a more realistic example in three dimensions we perform a simulation of a toroidal shell parametrized by 
\[
 x(\theta,\phi,r) = (R+r \cos(\theta)) \cos(\phi), \ \  y(\theta,\phi,r)= (R+r \cos(\theta)) \cos(\phi), \ \  z(\theta,\phi,r) = r \sin(\theta).  
 \]
Here we set $R = 4,$  and let the partial toroidal shell  occupy the volume $1\le r \le 2,$ and $0 \le \phi , \theta \le \pi$. The surfaces at $r=1$ and $r = 2$ are free, and we impose homogeneous Dirichlet conditions on all other boundaries.

We force the problem by 
\[
\bff =  \frac{\sqrt{\sigma^3}}{20}
 \left( \begin{array}{c}
1 \\
1 \\
1
\end{array}
 \right) e^{-\zeta^2} \cos(\omega t), \ \  \sigma = 100 \omega, \ \  \zeta^2 = 0.5 \sigma ((x-4)^2+(y-0.5)^2+(z-1)^2. 
\]

We consider three cases of increasing difficulty: 1) $\omega = 5.1234$ with a grid of $400 \times 120 \times 40$ points; 
2) $\omega =10.2468$ with a grid of $400 \times 120 \times 40$ points; 3) $\omega = 20.276$
with a grid consisting of $1600 \times 480 \times 160$ points. The largest computation thus solves the elastic Helmholtz problem for a system of equations with roughly $3.6\cdot 10^8$ degrees of freedom. We use the conjugate residual method, set a residual tolerance of $10^{-5}$ and obtain convergence in 1490, 2875, and 3713 iterations for cases 1), 2), and 3), respectively.

\begin{figure}[H]
\graphicspath{{figures/}}
\begin{center}
\includegraphics[width=0.89\textwidth,trim={0.0cm 1.65cm 0.0cm 1.7cm},clip]{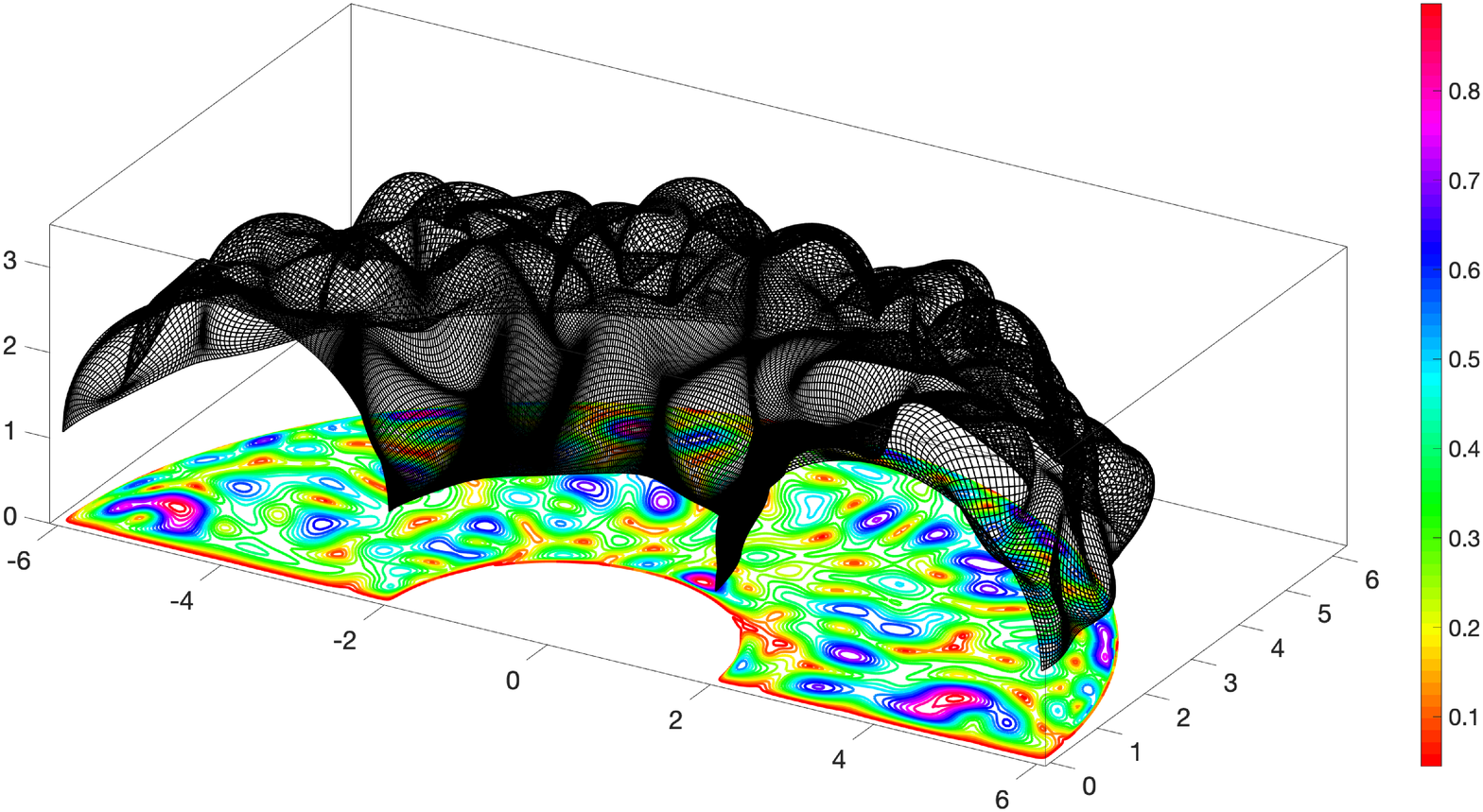}
\includegraphics[width=0.89\textwidth,trim={0.0cm 1.65cm 0.0cm 1.7cm},clip]{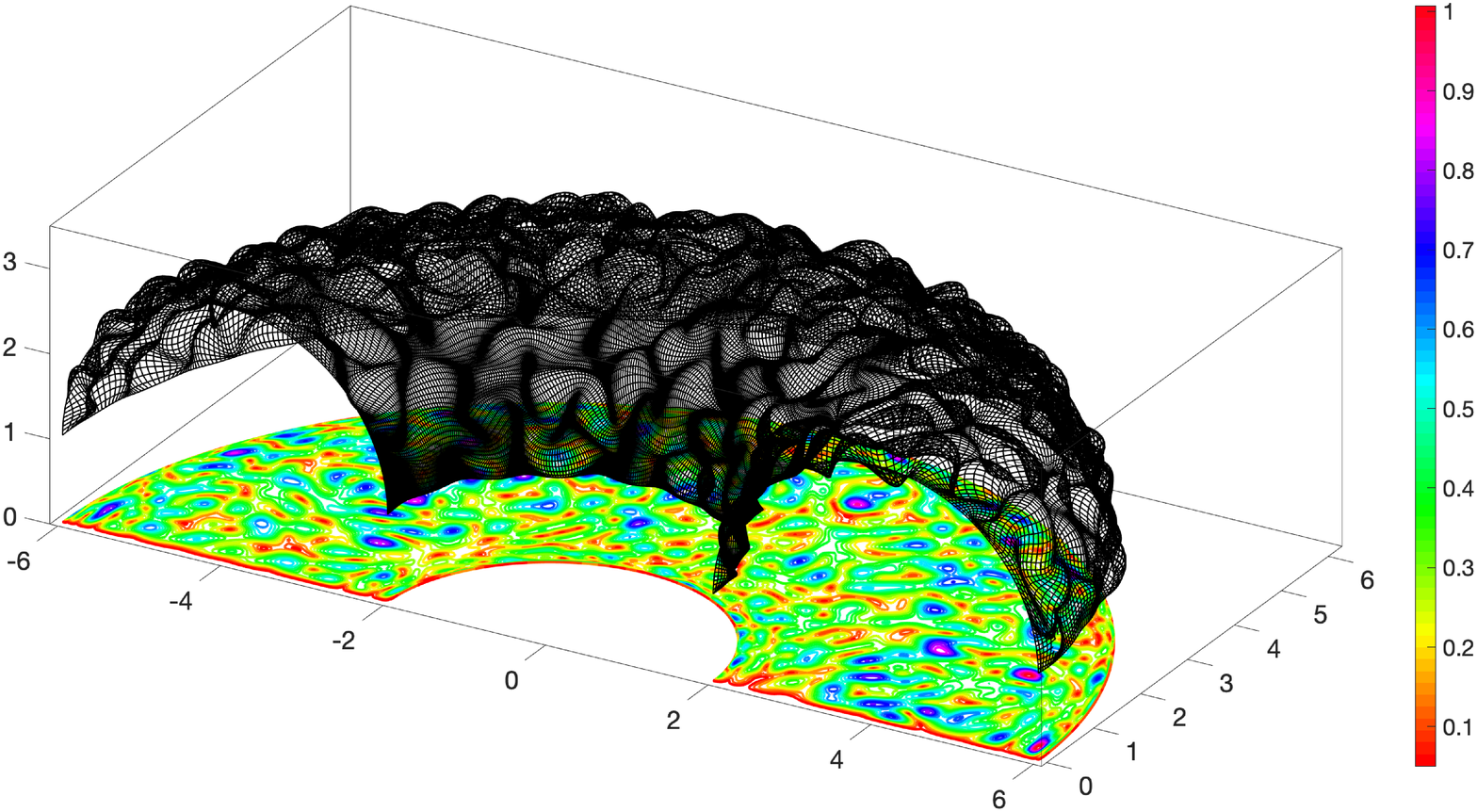}
\includegraphics[width=0.89\textwidth,trim={0.0cm 1.65cm 0.0cm 1.7cm},clip]{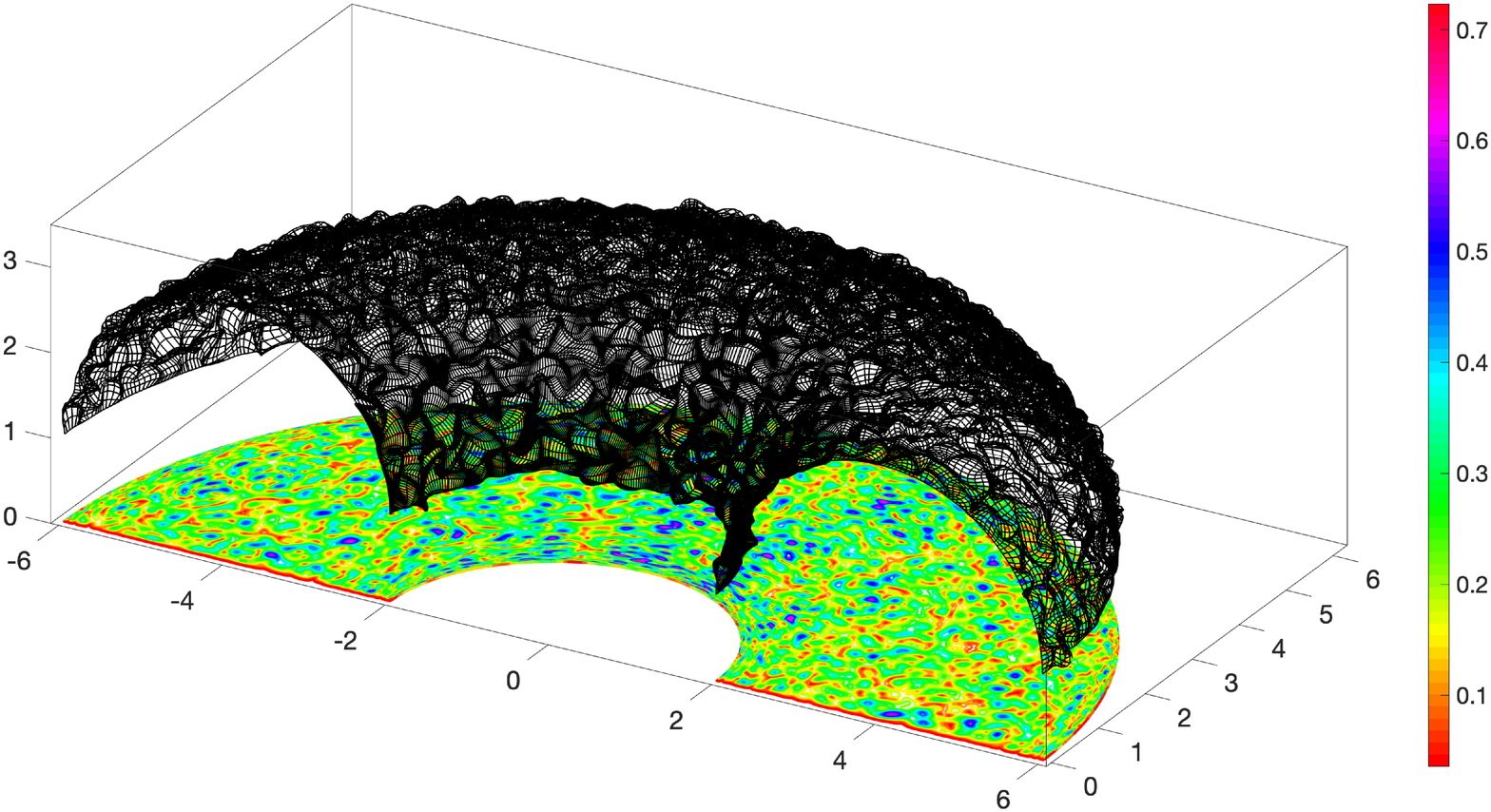}
\caption{The solution in the toroidal shell for (Left) $\omega = 5.1234$, (Middle) $\omega = 10.2468$, and (Right) $\omega = 20.276$. The projection onto the xy-plane is the magnitude of the displacement on the outermost free surface $r=2$. In black we display the (scaled) displaced mesh for $r=2$.
 \label{fig:torus}}
\end{center}
\end{figure}

The converged solutions are displayed in Figure \ref{fig:torus}. The projection onto the xy-plane is the magnitude of the displacement $\sqrt{u^2+v^2+w^2}$ on the outermost free surface, $r=2$. The mesh is the grid for that outermost surface with the (scaled) displacements added to the grid coordinates. 

\section{Conclusion}\label{sec::Conclusion}
In this paper we applied the WaveHoltz iteration, a time-domain 
Krylov accelerated fixed-point iteration, to the solution of the elastic Helmholtz equation
for interior problems with Dirichlet and/or free surface boundary
conditions. With symmetric discretizations, the iteration results 
in a positive definite and symmetric
matrix with a condition number that is bounded as the discretization parameter $h\to 0$, 
a notable advantage over direct discretizations
of the elastic Helmholtz equation which typically lead to highly indefinite
systems with condition numbers that grow as $h^{-2}$. In this work we have also proposed corrected
time-stepping schemes and demonstrated that their use in the WaveHoltz iteration
completely removes time discretization errors.

Furthermore, we have introduced a new implicit time-stepping
scheme for the El-WaveHoltz method, which
can offer some 
advantage over an explicit scheme, especially for high order, highly refined meshes with 
disparate element sizes. We believe that the implicit method could 
also be advantageous for anisotropic problems, and could potentially be 
an avenue for constructing polynomial/rational preconditioners due to the 
small number of time-steps afforded by an implicit scheme.

Finally, here we have only considered the energy conserving problem.
In the future, we will revisit the elastic 
Helmholtz problem with impedance/absorbing boundary conditions 
which are a hallmark of scattering and seismic applications.

\section*{Acknowledgements}
This work was supported in part by the NSF Grants DMS-1913076 and DGE-1650115; and in part 
by STINT initiation grant IB2019--8154. Any conclusions or recommendations expressed in this 
paper are those of the authors and do not necessarily reflect the views of the NSF.

\section*{Appendix}
\appendix

\section{Time-step Restriction}\label{appendix::RestrictionPlot}
To understand how restrictive the requirement \eqref{eqn::alpha} is, 
we plot $\alpha$ for various values of 
$\omega \Delta t$ below in Figure~\ref{fig::stability_plot}.
\begin{figure}[H]
\graphicspath{{figures/}}
\begin{center}
\includegraphics[width=0.48\textwidth]{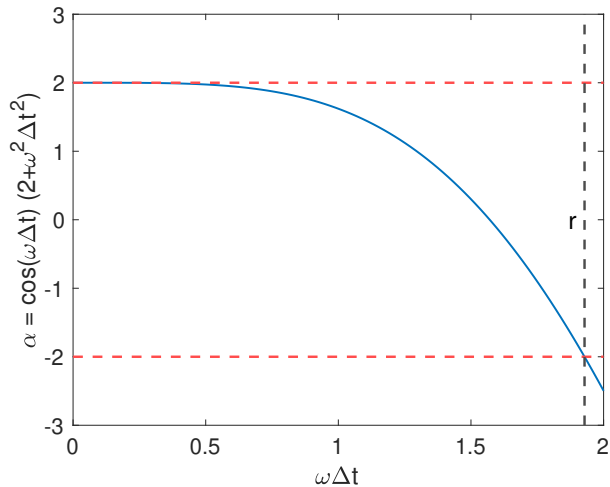}
\caption{Values of $\alpha = \cos(\omega \Delta t) (2 + \omega^2 \Delta t^2)$ for 
values of $\omega\Delta t$ in the interval $[0,2]$. The red lines 
indicate the desired bound on $\alpha$, and the black line indicates
the maximum allowable value of $\omega\Delta t$ at 
$r \approx 1.93$. \label{fig::stability_plot}}
\end{center}
\end{figure}
From Figure~\ref{fig::stability_plot} we see that $|\alpha| <2$ 
for $\Delta t < r/\omega$ where $r \approx 1.93$. 
This choice of the time-step corresponds to a requirement of at least four 
time-steps per iteration. However, the WaveHoltz kernel $2(\cos(\omega t_n)-1/4)/T$
evaluates to a constant if four time-steps are taken for the forward solve. Thus
at least five time-steps are needed for stability.

\section{Accuracy of the Discrete Filter Transfer Function}\label{appendix::TrapQuad}
\begin{figure}[H]
\graphicspath{{figures/}}
\begin{center}
\includegraphics[width=0.45\textwidth]{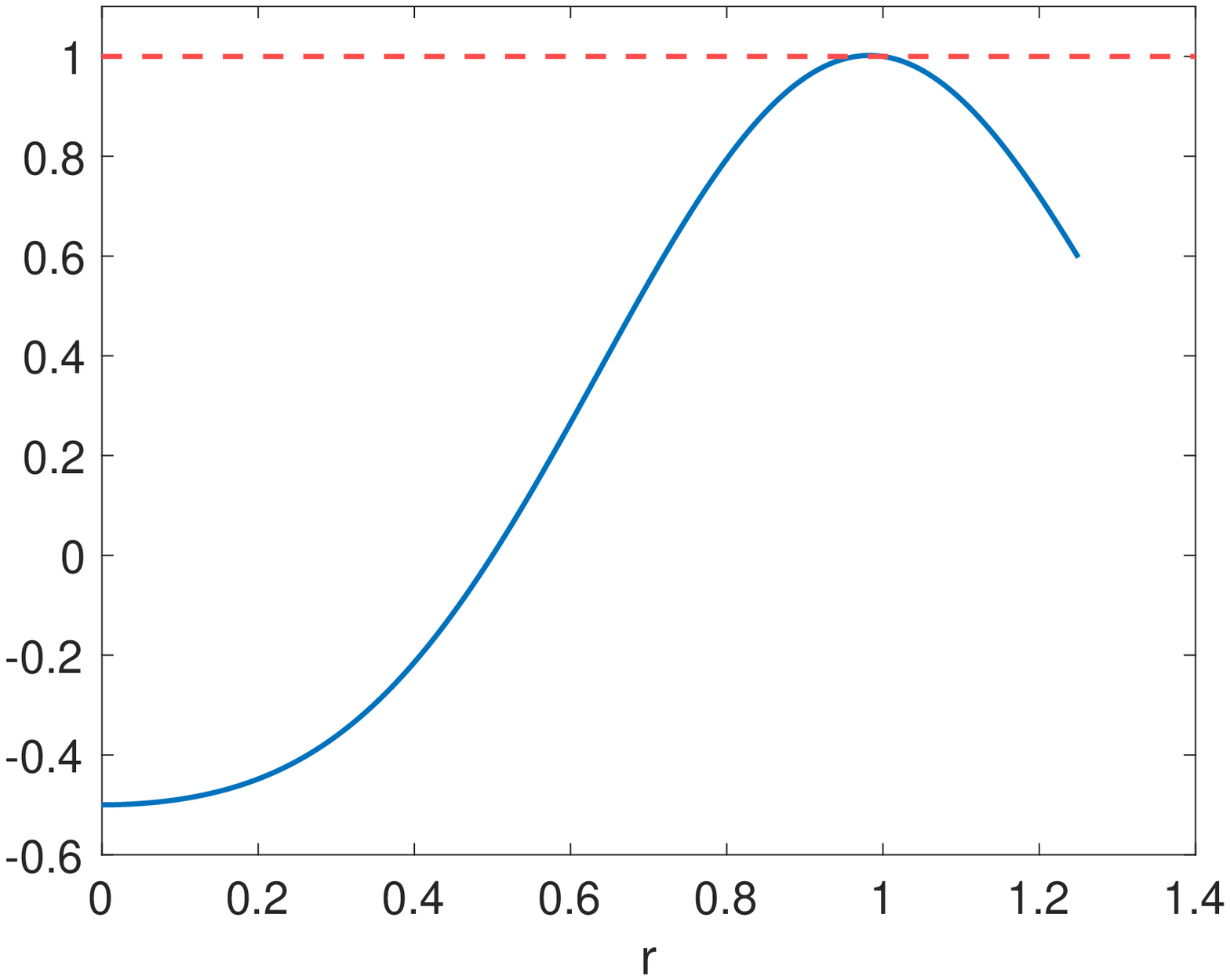}
\includegraphics[width=0.45\textwidth]{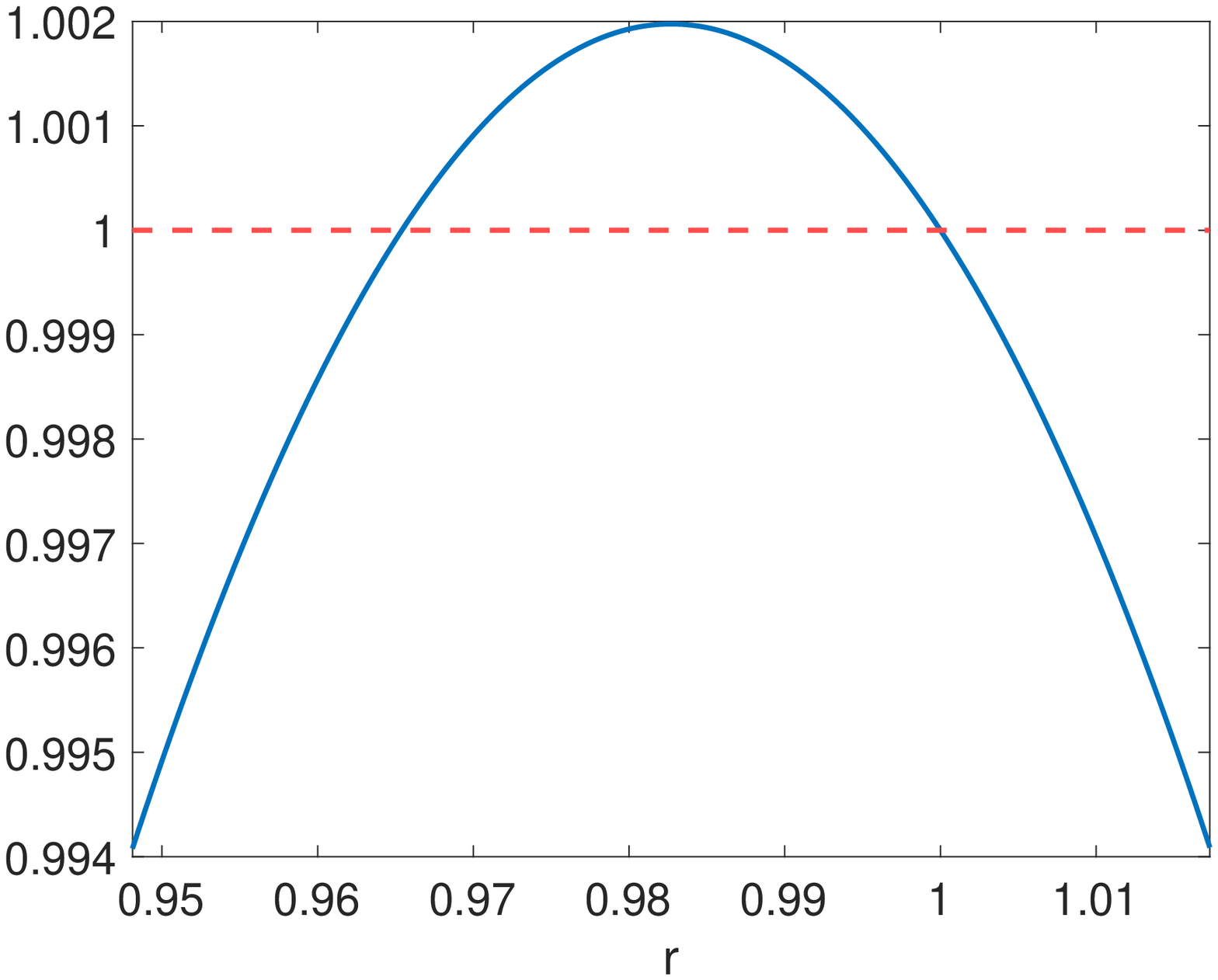}
\caption{A plot of the discrete filter function using five 
time-steps $0\le r\le 5/4$. On the left we plot the full range
of values of $r$, and on the right we zoom in close to resonance, i.e. $r=1$.
\label{fig::discrete_filter_func}}
\end{center}
\end{figure}

Consider the (continous) rescaled filter transfer function,
$$
  \bar{\beta}(r) := \beta(r\omega)
 = \frac{2}{T}\int_0^T\left(\cos(\omega t)-\frac14\right)\cos(r\omega t) dt = 
  \frac{1}{\pi}\int_0^{2\pi}\left(\cos(t)-\frac14\right)\cos(rt) dt,
$$
with discrete analogue (via trapezoidal rule)
$$
  \bar \beta_h(r) = 
   \frac{\Delta t}{\pi}
   \sum_{n=0}^M\eta_n 
   \cos(r t_n)\left(\cos(t_n)-\frac14\right),
\qquad \eta_n = \begin{cases}
   \frac12,& \text{$n=0$ or $n=M$},\\
   1, & 0<n<M.
   \end{cases}  
$$
Let us now take a look at the rescaled discrete filter function,
$\bar \beta_h(r)$. It is sufficient to consider only the range 
$r \in[0,k/4]$, which we plot in Figure~\ref{fig::discrete_filter_func}.

From Figure~\ref{fig::discrete_filter_func} we see that it is possible 
to integrate and get eigenvalues of the WHI operator to be larger than one
for a small range near resonance, $r=1$. To get a sense of the size of this gap,
we perform a simple bisection where we find the leftmost point $r^*<1$ such that
$\beta_h(r^*)=1$ for increasing number of quadrature points $k=5,6,\dots,100$.
We plot the result below in Figure~\ref{fig::discrete_filter_bound}.
\begin{figure}[H]
\graphicspath{{figures/}}
\begin{center}
\includegraphics[width=0.45\textwidth]{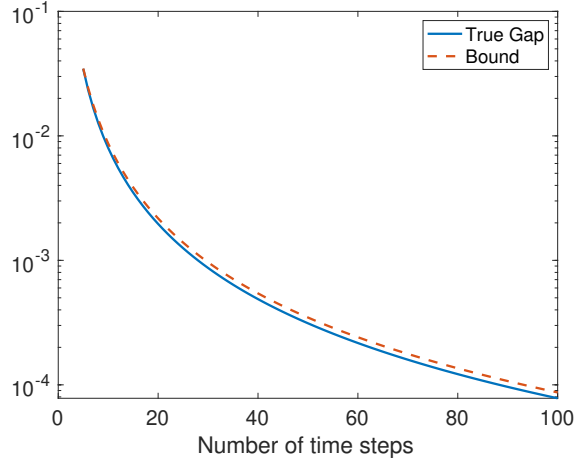}
\caption{A bound on the gap from resonance that creates problematic 
modes. The blue curve is the true gap, $1-r^*$, and the dotted
red curve is a proposed bound.
\label{fig::discrete_filter_bound}}
\end{center}
\end{figure}
From Figure~\ref{fig::discrete_filter_bound} we see that,
perhaps unsurprisingly, the gap shrinks with increasing number
of quadrature points. The curve in red in Figure~\ref{fig::discrete_filter_bound}
indicates the bound
  \begin{align*}
    1-r^* \le 0.022 \cdot \Delta t^2,
  \end{align*}
so that we see that this gap shrinks as $\Delta t^2$. Moreover, if 
$|1-r| \ge 0.022 \cdot \Delta t^2$ then $|\beta_h(r)| < 1$.
For $r = \tilde \lambda_j/\omega$, we may thus obtain convergence of the iteration if it 
can be guaranteed the time-step is 
chosen such that \mbox{$\tilde \lambda_j \not\in [\omega(1-0.022\cdot\Delta t^2),\omega]$}.

\bibliographystyle{siam}
\bibliography{garcia}
\end{document}